\documentclass[11pt]{article}

\usepackage{amsfonts}

\usepackage{CJK}
\usepackage{amsfonts,amssymb,amsmath,mathrsfs,multirow,booktabs}
\usepackage{xcolor,soul}
\usepackage{color,latexsym,amsfonts}
 \newcommand{\qed}{\hfill\rule{2mm}{3mm}\vspace{4mm}}

 \newtheorem{theorem}{Theorem}[section]
 \newtheorem{lemma}[theorem]{Lemma}
 \newtheorem{corollary}[theorem]{Corollary}
 \newtheorem{proposition}[theorem]{Proposition}
 \newtheorem{example}[theorem]{Example}
 \newtheorem{Definition}[theorem]{Definition}
 \newtheorem{remark}[theorem]{Remark}
 \newtheorem{condition}[theorem]{Condition}
 \newtheorem{conjecture}[theorem]{Conjecture}

 \def\blemma{\begin{lemma}\sl{}\def\elemma{\end{lemma}}}
 \def\btheorem{\begin{theorem}\sl{}\def\etheorem{\end{theorem}}}
 \def\bcorollary{\begin{corollary}\sl{}\def\ecorollary{\end{corollary}}}
 \def\bdefinition{\begin{Definition}\sl{}\def\edefinition{\end{Definition}}}
 \def\bproposition{\begin{proposition}\sl{}\def\eproposition{\end{proposition}}}
 \def\bremark{\begin{remark}\sl{}\def\eremark{\end{remark}}}

 \def\beqlb{\begin{eqnarray}}\def\eeqlb{\end{eqnarray}}
 \def\beqnn{\begin{eqnarray*}}\def\eeqnn{\end{eqnarray*}}

 \def\mbb{\mathbb}\def\mbf{\mathbf}

 \def\<{\langle}\def\>{\rangle}

 \def\ar{&\!\!}

 \def\eqref#1{{\rm(\ref{#1})}}

 \def\proof{\noindent{\it
 Proof.~}}\def\qed{\hfill$\Box$\medskip}

\def\e{{\mbox{\rm e}}}
\def\<{\left<}\def\>{\right>}

 \def\mbb{\mathbb}
  \def\mbf{\mathbf}
\newcommand{\dd}{\mathrm{d}}

\newfam\msbmfam\font\tenmsbm=msbm10\textfont
\msbmfam=\tenmsbm\font\sevenmsbm=msbm7
\scriptfont\msbmfam=\sevenmsbm

\def\<{\left<}\def\>{\right>}
\def\({\left(}\def\){\right)}

\begin{document}

\centerline{\Large\bf Extinction behaviour for }
\medskip
\centerline{\Large\bf competing continuous-state population dynamics}

\bigskip

\centerline{
Jie Xiong\footnote{Department of Mathematics and SUSTech International center for Mathematics, Southern University of Science \& Technology, Shenzhen, China.
Supported by National Key R\&D Program of China (No.~2022YFA1006102).  Email: xiongj@sustech.edu.cn},
Xu Yang\footnote{School of Mathematics and Information Science,
North Minzu University, Yinchuan, China. Supported by
NSFC (No.~12471135) and NSF of Ningxia (No.~2024AAC05056). Email: xuyang@mail.bnu.edu.cn. Corresponding author.} and Xiaowen Zhou\footnote{Department of Mathematics and Statistics, Concordia
University, Montreal, Canada.
Supported by NSERC (RGPIN-2021-04100).
Email: xiaowen.zhou@concordia.ca.}}

\bigskip\bigskip

{\narrower{\narrower

\noindent{\bf Abstract}
We consider a system of two stochastic differential equations (SDEs) with  competing two-way interactions driven by Brownian motions  and spectrally positive $\alpha$-stable random measures. Such a SDE system can be identified as a Lotka-Volterra type population model.
We find nearly sharp conditions for one of the population to become extinct or extinguished.

\bigskip


\bigskip

\par}\par}

\section{Introduction and main result}\label{intro}

\setcounter{equation}{0}

Extinction behavior is a key topic in the study of population models. For the continuous-state branching process (CSBP for short) arising as scaling limits of the Galton-Watson branching processes, a sufficient and necessary condition, called Grey's condition, is obtained  in \cite{Grey74}. Generalizations of CSBPs have been introduced since then. In \cite{LiP19} a class of CSBPs with state-dependent branching mechanism  was obtained by a Lamperti type time change of spectrally positive L\'evy processes stopped at the first time of hitting $0$. Integral tests for the extinction of such nonlinear CSBPs are obtained in \cite{LiP19} and \cite{LiZh21}.
The above nonlinear CSBPs were further generalized in \cite{LYZh} to solutions to stochastic differential equations (SDEs) driven by Brownian motion and spectrally positive Poisson random measures. Rather sharp conditions are obtained on extinction/non-extinction for these nonlinear CSBPs, see also \cite{MYZ21} for results on critical cases and  \cite{baiyang} for the nonlinear CSBPs with Neveu's branching. Similarly to extinction, a process becomes extinguished if it converges to $0$ but never reaches $0$. A criterion in terms of integral test is found in \cite{LTZ} for the nonlinear CSBP as time changed spectrally positive L\'evy process to be extinguished.

The study of extinction behaviour belongs to the class of problems of boundary classification for nonnegative valued Markov processes where the concern is whether $0$ or $\infty$ is an entrance, exit, or neutral boundary.

The main approach to developing the criteria for extinction/non-extinction
on the above work is an adaption of the approach for  Chen's criteria on the uniqueness problem of Markov jump processes. These Chen's criteria are first proposed in \cite{Chen1986a,Chen1986b} and can also be found in
\cite[Theorems 2.25 and 2.27]{Chen04}.
Such an approach typically involves identifying an appropriate test function that is applied to the infinitesimal generator of the nonlinear CSBP and proving the desired result using a martingale argument. It also finds successful applications in studying other boundary behaviours for Markov processes such as explosion/nonexplosion and coming down from infinity versus staying infinite. We refer to \cite{MT93} for a similar approach to studying boundary behaviours for Markov processes.

Compared with one-dimensional models, the study of two-dimensional interacting population dynamics turns out to be much more challenging, and there are  relatively few results available. We refer to \cite{DuDangYin} and \cite{BaoShao}  for related permanence/extinction results on stochastic predator-prey models as coupled diffusions and to \cite{Bao} for a multi-dimensional interacting Lotka-Volterra system subject to a common one-dimensional jump diffusion type environmental noise.

In \cite{RXYZ19} a stochastic Lotka-Volterra type population dynamical system $(X,Y)$ was proposed as a solution to two-dimensional SDEs with one sided interaction where the dynamics of process $Y$ is affected by $X$. When the continuous-state population $X$ is extinguished, quite sharp conditions are found for the population $Y$ to become extinct or extinguished.   The main approach of \cite{RXYZ19} is again an adaptation of that of Chen's criteria  technique from one-dimensional processes to two-dimensional processes.

Given the previous result on two-dimensional SDE  models with one-sided interaction, it is natural to consider models with two-way interactions. We are not aware of systematic work on boundary behaviours of two-dimensional stochastic processes, not even two-dimensional interacting diffusions. As a first attempt, in this paper we start with the extinction behavior of the following SDE system:
\begin{equation}\label{a1.1}
	\left\{
	\begin{aligned}		
X_t &=X_0-\int_0^t[a_1X_s^{p_1+1}+\eta_1X_s^{\theta_1}Y_s^{\kappa_1}] \dd s+\int_0^t\sqrt{2a_2 X_s^{p_2+2}}\dd B_1(s)\\
&~~~~
+\int_0^t\int_0^\infty \int_0^{a_3X_{s-}^{p_3+\alpha_1}}z\tilde{N}_1(\dd s,\dd z,\dd u),\\
		Y_t &=Y_0-\int_0^t[b_1 Y_s^{q_1+1}+\eta_2Y_s^{\theta_2}X_s^{\kappa_2}]\dd s+\int_0^t\sqrt{2b_2 Y_s^{q_2+2}}\dd B_2(s)\\
&~~~~
+\int_0^t\int_0^\infty \int_0^{b_3 Y_{s-}^{q_3+\alpha_2}}z\tilde{N}_2(\dd s,\dd z,\dd u),
	\end{aligned}
	\right.
\end{equation}
where the constants $X_0,Y_0,\kappa_i,\eta_i>0$ and
$\theta_i,a_j,b_j,p_j,q_j\ge0$ for $i=1,2$ and $j=1,2,3$.
For $i=1,2$, $(B_i(t))_{t\ge0}$ is a Brownian motion and
$\tilde{N}_i(\dd s,\dd z,\dd u)$ is a compensated Poisson
random measure with intensity $\dd s\mu_i(\dd z)\dd u$.
Here
$\mu_i(\dd z)=\frac{\alpha_i(\alpha_i-1)}{\Gamma(\alpha_i)\Gamma(2-\alpha_i)}z^{-1-\alpha_i}\dd z,z>0$,
for $\alpha_i\in(1,2), i=1,2,$ and  $\Gamma$ denotes the Gamma function.
We always assume that $a_2+a_3>0$ and $b_2+b_3>0$.
We also assume that
$(B_1(t))_{t\ge0}$, $(B_2(t))_{t\ge0}$, $\{\tilde{N}_1(\dd s,\dd z,\dd u)\}$, and $\{\tilde{N}_2(\dd s,\dd z,\dd u)\}$ are independent of each other.

The extinction-extinguishing dichotomy is studied in \cite{RXYZ19}
 for the SDE system \eqref{a1.1} satisfying $\eta_1=0$.
If $p_1=q_1=p_2=q_2=0$, $\theta_1=\theta_2=\kappa_1=\kappa_2=1$,
$a_3=b_3=0$, and $a_1,b_1<0$, then
SDE \eqref{a1.1} is also called the competitive Lotka-Volterra model in random environments in \cite{Bao} and the references, and \cite{Evans}
considers the case where the driven noise $((B_1(t),B_2(t)))_{t\ge0}$ is a correlated two-dimensional Brownian motion.
In \cite{Cattiaux}, the quasi-stationary distribution is studied
for $p_1=q_1=0$, $p_2=q_2=-1$, $\theta_1=\theta_2=\kappa_1=\kappa_2=1$,
$a_3=b_3=0$, and $a_1,b_1<0$ in \eqref{a1.1}.

In the SDE system \eqref{a1.1},
for general measures $\mu_i$ satisfying $\int_0^\infty (z\wedge z^2)\mu_i(\dd z)$, and
with $a_1x^{p_1+1}$, $2a_2 x^{p_2+2}$, $a_3x^{p_3+\alpha_1}$,
$b_1y^{q_1+1}$, $2b_2y^{q_2+2}$, $b_3y^{q_3+\alpha_2}$, $\eta_1X_s^{\theta_1}Y_s^{\kappa_1}$ and $\eta_2Y_s^{\theta_2}X_s^{\kappa_2}$
replaced
by nonnegative functions $p_1(x)$, $p_2(x)$, $p_3(x)$,
$q_1(y)$, $q_2(y)$, $q_3(y)$, $\eta_1(x,y)$ and $\eta_2(x,y)$, respectively,
we can also study the extinction behaviour by the same technique in this paper.
For simplicity and readability, we only consider this special form of \eqref{a1.1}.

For process $Z=(Z_t)_{t\ge0}$ and $u>0$ define stopping times
$$\tau_u^-(Z)=\inf\{t\ge0:Z_t \le u\},\,\, \tau_u^+(Z)=\inf\{t\ge0:Z_t \ge u\},$$
$$\tau_0^-(Z):=\inf\{t\ge0:Z_t=0\}\,\,\, \text{and}\,\,\, \tau_\infty^+(Z):=\inf\{t\ge0:Z_t=\infty\}$$
 with the convention $\inf\emptyset=\infty$.

We first present the definition of solution to SDE \eqref{a1.1}, which is
defined up to the time when any of the two processes first hits $0$ or reaches infinity. For the solution $((X_t,Y_t))_{t\ge0}$ to SDE \eqref{a1.1}
let $X=(X_t)_{t\ge0}$ and $Y=(Y_t)_{t\ge0}$.

\bdefinition\label{def}
A two-dimensional c\`adl\`ag
process $((X_t,Y_t))_{t\ge0}$ is a solution to SDE \eqref{a1.1} if
it is defined in the same filtered probability space as that for Brownian motions $(B_1(t))_{t\ge0}$
and $(B_2(t))_{t\ge0}$ and compensated Poisson
random measures $\{\tilde{N}_1(\dd s,\dd z,\dd u)\}$ and $\{\tilde{N}_2(\dd s,\dd z,\dd u)\}$, and satisfies
SDE \eqref{a1.1} up to $\sigma_n:=\tau_{1/n}^-(X)\wedge\tau_{1/n}^-(Y)
\wedge\tau_{n}^+(X)\wedge\tau_{n}^+(Y)$ for each $n\ge1$ and
$X_t=\limsup_{n\to\infty}X_{\sigma_n-}$
and $Y_t=\limsup_{n\to\infty}Y_{\sigma_n-}$
for $t\ge\lim_{n\to\infty}\sigma_n$.
\edefinition

The definition \ref{def} implies that
$0$ is an absorbing state and the process $((X_t,Y_t))_{t\ge0}$ stays in
$(\frac1n,\infty)\times(\frac1n,\infty)$ before $\sigma_n$
since $X$ and $Y$ have upwards jumps.
Thus, the solution is nonnegative and
the definition \ref{def} of the solution to SDE \eqref{a1.1}
allows weaker conditions for the uniqueness of the solution. In particular, the existence and pathwise uniqueness are valid following the same arguments as in \cite[Lemma A.1]{RXYZ19}.
Throughout this paper, we assume that the  c\`adl\`ag process $((X_t,Y_t))_{t\ge0}$ is a unique solution to (1.1), and consequently the process $((X_t,Y_t))_{t\ge0}$ has the strong Markov property.  We
also assume that $X_0,Y_0>0$ are deterministic and that all stochastic processes are defined on the same filtered probability space $(\Omega,\mathscr{F},\mathscr{F}_t,\mathbf{P})$.
We use $\mathbf{E}$ to denote the corresponding expectation.

The negatively interacting term represents the negative effects of competition between the two populations. If $\eta_1=0$ in (\ref{a1.1}), then process $X$ is a CSBP with a nonlinear branching mechanism. If $p_1=0, p_2=-1$ and $p_3=1-\alpha_1$, then $X$ further reduces to a CSBP that is a scaling limit of the classical Galton-Watson branching process.  We refer to \cite{Ky14} and \cite{Li} for a complete introduction to CSBPs. Note that the values of the parameters are chosen so that when there is no interaction, i.e. $\eta_i=0$ for $i=1,2$, neither $X$ nor $Y$ can reach $0$ with a positive probability (see \cite[p.2535]{LYZh}). Therefore,  we would not consider CSBPs for our purpose.

The processes $X$ and $Y$ become extinguished when $\eta_1=\eta_2=0$
(see \cite[p.2535]{LYZh}).
The goal of this paper is to find as sharp as possible conditions
on the interaction terms for the extinction and extinguishing behaviours
of the mutually interacting nonlinear CSBP \eqref{a1.1}.
Write
$\tau_0^-:=\tau_0^-(X)\wedge\tau_0^-(Y)$.
Let
 \beqnn
p:=\min\big\{p_i, \, i\in \{j: a_j\neq 0, j=1,2,3\}\big\}, \quad
q
:=
\min\big\{q_i, \,i\in \{j: b_j\neq 0, j=1,2,3\}\big\},
 \eeqnn
i.e., $p$ is the minimum of those $p_i$ with $a_i\neq 0$ and $q$ is the minimum of those $q_i$ with $b_i\neq 0$.
Let
 \beqnn
a:=
a_11_{\{p_1=p\}}+a_21_{\{p_2=p\}}
+a_3 1_{\{p_3=p\}},
\quad
b:=
b_11_{\{q_1=q\}}+b_21_{\{q_2=q\}}
+b_3 1_{\{q_3=q\}}.
 \eeqnn

The following theorems are the main results of the paper.
For the nonlinear one-dimensional CSBP (with $\eta_1X_s^{\theta_1}Y_s^{\kappa_1}$
replaced by $\eta_1X_s^{\theta_1}$ in \eqref{a1.1}),
the process cannot reach $0$ iff $\theta_1\ge1$ (see \cite[p.2535]{LYZh}).
The first theorem  shows that this conclusion still holds
for the two-dimensional model \eqref{a1.1}, and $Y_s^{\kappa_1}$
in the interaction term
$\eta_1X_s^{\theta_1}Y_s^{\kappa_1}$ does not affect the extinction of $X$
when $\theta_1\ge1$.

\btheorem\label{t0.01}
$\mbf{P}\{\tau_0^-<\infty\}=0$
iff $\theta_1,\theta_2\ge1$.
Moreover,
$\mbf{P}\{\tau_0^-(X)<\infty\}=0$
for $\theta_1\ge1$
and
$\mbf{P}\{\tau_0^-(Y)<\infty\}=0$
for $\theta_2\ge1$.
\etheorem

For the one-dimensional nonlinear CSBP (with $\eta_1X_s^{\theta_1}Y_s^{\kappa_1}$
replaced by $\eta_1X_s^{\theta_1}$ in \eqref{a1.1}),
the process reaches $0$ with probability one when $\theta_1<1$ (see
\cite[Corollary 2.5(ii) and Proposition 2.6(i)]{LYZh}).
However, this conclusion does not hold for the two-dimensional model \eqref{a1.1}.
The following theorems state the sharp conditions on extinction with probability one
and extinction with probability strictly between zero and one.
The following theorem is stated for $Y$, and similar results hold for $X$ by symmetry.

\btheorem\label{t0.02}
Assuming that $\theta_1\ge1$ and $0\le\theta_2<1$,
if one of the following holds:
\begin{itemize}
\item[{\normalfont(a)}]
$\theta_1<1+\frac{\kappa_2(q-\kappa_1)}{q+1-\theta_2}$;
\item[{\normalfont(b)}]
$p<\frac{\kappa_2q}{q+1-\theta_2}$;
\item[{\normalfont(c)}]
$p=q=0$ and $b/a<\kappa_2/(1-\theta_2)$;
\item[{\normalfont(d)}]
$\theta_1=1$, $q=\kappa_1$, and $b/\eta_1<\kappa_2/(\kappa_1+1-\theta_2)$,
\end{itemize}
then,
$0<\mbf{P}\{\tau_0^-(Y)<\infty\}<1$.
\etheorem

We note that conditions (a) and (b) concern the noncritical case in which
the result is only affected by the powers of the terms involved, while conditions (c) and (d) concern the critical case in which the values of coefficients $a_i, b_i$ and $\eta_i$ can affect the conclusion.


\btheorem\label{t0.03}
Assume that $\theta_1\ge1$ and $0\le\theta_2<1$ and
\begin{itemize}
\item[{\normalfont(a)}]
$\theta_1>1+\frac{\kappa_2(q-\kappa_1)}{q+1-\theta_2}$,
\end{itemize}
if one of the following further holds:
\begin{itemize}
\item[{\normalfont(b)}]
$p>\frac{\kappa_2q}{q+1-\theta_2}$;
\item[{\normalfont(c)}]
$p=q=0$ and $b/a>\kappa_2/(1-\theta_2)$,
\end{itemize}
then, $\mbf{P}\{\tau_0^-(Y)<\infty\}=1$.
\etheorem

\bremark
(i) Note that neither $X$ nor $Y$ in (\ref{a1.1}) can reach $0$ if $\eta_1=\eta_2=0$ by \cite[Theorem 2.3(i)]{LYZh}. Therefore, extinction is caused by the interaction between processes $X$ and $Y$.
If $\eta_1=0$, the results of Theorems \ref{t0.01}--\ref{t0.03}
are given in \cite[Example 1.12]{RXYZ19}.

(ii) The extinction behaviours for the critical cases in the last two theorems remain unknown, and we leave them as an {\em open} problem.

(iii)
Under the setup of SDE (\ref{a1.1}), both processes $X$ and $Y$ converge to $0$, which facilitates the proofs for extinction/extinguishing. It would be more challenging to show similar results when both of the interacting terms are positive. The approach of this paper also finds successful applications in showing the criteria for explosion/non-explosion of the solution to one-dimensional SDE; see \cite{MYZ21} and \cite{MZ24}. We expect that it can also be modified to show similar criteria for two-dimensional models with interaction.
\eremark

Theorem \ref{t0.01} is proved in Lemmas \ref{t3.3} and \ref{t3.4}, respectively.
The proof for Lemma \ref{t3.3} is given by using a non-extinction criterion established in Proposition \ref{t3.1b}.
The proof of assertion ``$\mbf{P}\{\tau_0^-(Y)<\infty\}>0$"
in Theorem \ref{t0.02} is given in Lemmas \ref{t3.3} and \ref{t3.4}.
The proof of Lemma \ref{t3.4} is essentially the same
as that of \cite[Theorem 1.5]{RXYZ19} using an extinction criterion established
in \cite{RXYZ19} (see Proposition \ref{t0.3})
and a complicated exponential function (see \eqref{2.0}).
The proof of
Theorem \ref{t0.03} is established by applying an extinction criterion
established
in \cite{RXYZ19} (see Corollary \ref{cor2.1}) to an exponential test
function on a ratio of the two processes.

The assertion in Theorem \ref{t0.02} that the probability of extinction is strictly less than one is proved by a new method different from that of
\cite{RXYZ19}.
To this end, we first prove a criterion for the extinction  probability to be strictly less than one
for small initial values by considering a ratio $U$ of processes $Y$ and $X^\beta$ for which both $X_0$ and $Y_0$ are small but $U_0$ is large; see Proposition \ref{Pro2.1}.
We then construct a test function and verify the criterion; see Subsection \ref{subsection3.3}. We also establish an irreducibility criterion in Proposition \ref{pro2.2a}, comparison theorem in Proposition \ref{t2.1} and apply the criterion   to show the result for any initial values; see
the end of the proof of Theorem \ref{t0.02}.

Throughout the paper
$C^2((0,\infty))$ and $C^2((0,\infty)\times(0,\infty))$
denote the spaces of second-order continuous differentiable functions  on $(0,\infty)$
and $(0,\infty)\times(0,\infty)$, respectively.
The remainder of the paper is arranged as follows. In Section \ref{Criteria},
we prove some preparatory  results on the finiteness of the hitting time $\tau^-_0$ and the irreducibility of a generic two-dimensional
process.  Section \ref{proof of Theorems} contains proofs for the main results.
Section \ref{Appendix}  contains two key lemmas for the proofs in Section \ref{proof of Theorems} and the proof for Proposition \ref{t2.1}.

\section{Criteria for the extinction behaviour}\label{Criteria}
\setcounter{equation}{0}

In this section, we establish some preliminary results for Theorems \ref{t0.01}--\ref{t0.03}.
In the following, let $((x_t,y_t))_{t\ge0}$ be a generic two-dimensional
process where $x=(x_t)_{t\ge0}$ and $y=(y_t)_{t\ge0}$ are two nonnegative
processes defined up to the minimum of their first times of hitting
zero or explosion. We assume that $x_0,y_0>0$ are deterministic,
and $t\mapsto x_t$ and $t\mapsto y_t$ are right continuous with left limits.

Let $\mathcal{L}$ denote the operator so that for each $g\in C^2((0,\infty)\times(0,\infty))$,
the process
 \beqlb\label{3.1}
t\mapsto M^g_{t\wedge\tau_{m,n}} \quad \mbox{is a local martingale},
 \eeqlb
where
 \beqnn
M^g_t:=g(x_t,y_t)-g(x_0,y_0)-\int_0^t\mathcal{L}g(x_s,y_s)\dd s
 \eeqnn
and $\tau_{m,n}:=\tau^-_{1/m}\wedge\tau^+_n$
with $\tau_{1/m}^-:=\tau^-_{1/m}(x)\wedge\tau^-_{1/m}(y)$
and $\tau_n^+:=\tau^+_{n}(x)\wedge\tau^+_{n}(y)$.
For SDE, the operator $\mathcal{L}$ can be obtained using It\^o's formula.
In this section, we assume that $\tau_0^-:=\tau_0^-(x)\wedge\tau_0^-(y)$.

The proofs for the propositions in this section are based on
\eqref{3.1} and martingale techniques.
We first state the non-extinction and extinction criteria, which generalize Chen's criteria for the uniqueness problem of Markov jump processes.
In applying such a criterion, the key is to identify an
appropriate test function that leads to the best possible result, but it
may not always allow for an obvious intuitive interpretation.

The following criterion is on non-extinction of process $((x_t,y_t))_{t\ge0}$
and the test function $g_n$ is often selected to be of logarithm type
or power function with negative power.

\bproposition\label{t3.1b}
Assuming that $\sup_{t\ge0} (x_t+y_t)<\infty$ a.s.,
 \beqlb\label{3.0}
\mbf{P}\{\tau_0^-(x)=\tau_0^-(y)<\infty\}=0,
 \eeqlb
and that for any  $n\ge1$,
there is a nonnegative function $g_n\in C^2((0,\infty)\times(0,\infty))$
and a constant $d_n>0$ such that
\begin{itemize}
\item[{\normalfont(i)}]
$\lim_{u\to0}g_n(u,v)=+\infty$ for all $v>0$,
\item[{\normalfont(ii)}]
$\mathcal{L}g_n(u,v)\le d_ng_n(u,v)$
for all $0<u,v\le n$,
\end{itemize}
then we have $\mbf{P}\{\tau_0^-(x)=\infty\}=1$.
\eproposition
\proof
By \eqref{3.1},
there are stopping times $\gamma_k$ so that $\gamma_k\to\infty$ as $k\to\infty$
almost surely
and $t\mapsto M^g_{t\wedge\tau_{m,n,k}}$
is a martingale, where $\tau_{m,n,k}:=\tau_{m,n}\wedge\gamma_k$.
Thus, for all $m,n,k\ge1$,
 \beqnn
\mbf{E}\big[g_n(x_{t\wedge \tau_{m,n,k}},y_{t\wedge\tau_{m,n,k}})\big]
 \ar=\ar
g_n(x_0,y_0)+\mbf{E}\Big[\int_0^{t\wedge\tau_{m,n,k}}
\mathcal{L}g_n(x_s,y_s)\dd s\big] \cr
 \ar=\ar
g_n(x_0,y_0)+\mbf{E}\Big[\int_0^{t}
\mathcal{L}g_n(x_{s},y_{s})
1_{\{s<\tau_{m,n,k}\}}\dd s\Big] \cr
 \ar\le\ar
g_n(x_0,y_0)+d_n\int_0^t
\mbf{E}\big[g_n(x_{s\wedge\tau_{m,n,k}},y_{s\wedge\tau_{m,n,k}})\big]\dd s,
 \eeqnn
where condition (ii) is used in the last inequality.
It follows from Gronwall's inequality that
 \beqnn
\mbf{E}\big[
g_n(x_{t\wedge\tau_{m,n,k}},y_{t\wedge\tau_{m,n,k}})\big]
\le g_n(x_0,y_0)\e^{d_nt}.
 \eeqnn
For convenience, we assume that $g(0,v):=\lim_{u\to0}g(u,v)$
and $g(v,0):=\lim_{u\to0}g(v,u)$ for $v>0$.
By Fatou's lemma,
 \beqlb\label{2.2}
 \ar\ar
\mbf{E}\Big[g_n(x_{t\wedge\tau_0^-(x)\wedge\tau_0^-(y)\wedge\tau_n^+},
y_{t\wedge\tau_0^-(x)\wedge\tau_0^-(y)\wedge\tau_n^+})\Big]=
\mbf{E}\Big[\liminf_{m,k\to\infty}
g_n(x_{t\wedge\tau_{m,n,k}},y_{t\wedge\tau_{m,n,k}})\Big] \cr
 \ar\ar\qquad\le
\liminf_{m,k\to\infty}\mbf{E}\big[
g_n(x_{t\wedge\tau_{m,n,k}},y_{t\wedge\tau_{m,n,k}})\big]
\le g_n(x_0,y_0)\e^{d_nt}.
 \eeqlb
Let $S_t^n:=\{\tau_0^-(x)<t\wedge\tau_0^-(y)\wedge\tau_n^+\}$.
If $\mbf{P}\{S_t^n\}>0$, then combining \eqref{2.2}, condition (i) and the fact $g_n\ge0$,
we have $g(0,v)=\infty$ for $v>0$ and
 \beqnn
\infty=\mbf{E}\Big[g_n(x_{\tau_0^-(x)},
y_{\tau_0^-(x)})
1_{S_t^n}\Big]\le g_n(x_0,y_0)\e^{d_nt}<\infty,
 \eeqnn
which is a contradiction.
Thus
 \beqnn
\mbf{P}\{\tau_0^-(x)\ge t\wedge\tau_0^-(y)\wedge\tau_n^+\}=1,\qquad t>0.
 \eeqnn
By the assumption, $\tau_n^+\to\infty$ as $n\to\infty$ and then
letting $n,t\to\infty$ we obtain $\tau_0^-(x)\ge\tau_0^-(y)$ almost surely.
Then one can reach the desired result by \eqref{3.0} and the definition of
$((x_t,y_t))_{t\ge0}$ before the minimum of their first times of hitting zero.
\qed

The following two criteria on the extinction of the process $((x_t,y_t))_{t\ge0}$
can be found in \cite{RXYZ19} and the key function $g$ is often the composition of a power function and an exponential function.
\bproposition\label{t0.3}
(\cite[Proposition 2.2]{RXYZ19})
Suppose that $\sup_{t\ge0} (x_t + y_t)<\infty$ a.s.
and in addition there exist a nonnegative function $g\in C^2((0,\infty)\times(0,\infty))$ and constants $d_n>0$ satisfying the following conditions:
\begin{itemize}
\item[{\normalfont(i)}] $0<\sup_{u,v>0} g(u,v)<\infty$;
\item[{\normalfont(ii)}] $\mathcal{L}g(u,v)\geq d_ng(u,v)$ for all
$0<u,v\le n$.
\end{itemize}
Then
$\mbf{P}\{\tau_0^-<\infty\}\ge g(x_0,y_0)/ \sup_{x,y>0}g(x,y)$.
\eproposition

\bcorollary\label{cor2.1}
(\cite[Corollary 2.3]{RXYZ19})
Suppose that $\sup_{t\ge0} (x_t+y_t)<\infty$ a.s.
and $g\in C^2((0,\infty)\times(0,\infty))$ is a nonnegative function
with $0<\sup_{u,v>0}g(u,v)<\infty$.
If there are constants $c_0,\varepsilon>0$ such that $x_0,y_0<\varepsilon$,
 \beqlb\label{2.4a}
\mathcal{L}g(u,v)\geq c_0g(u,v),\qquad 0<u,v<\varepsilon,
 \eeqlb
then
 \beqnn
\mbf{P}\{\tau_0^-\wedge\tau_\varepsilon^+<\infty\}
\ge g(u_0,v_0)/\sup_{u,v>0}g(u,v).
  \eeqnn
\ecorollary

The next result concerns the irreducibility of the process $(x,y)$,
i.e., it can  reach  any open subset of $(0,\infty)\times(0,\infty)$ with positive probability when the system starts from any  point in $(0,\infty)\times(0,\infty)$,
which is needed for the proof of Theorem \ref{t0.02}.

\bproposition\label{pro2.2a}
For $u_2>u_1>0$ and $v_2>v_1>0$ let $D:=[u_1, u_2]\times [v_1,v_2]$.
Let $0<\varepsilon_1\le u_1,0<\varepsilon_2\le v_1,\delta_1\ge u_2,\delta_2\ge v_2$, and $\varepsilon_1<x_0<\delta_1,\varepsilon_2<y_0<\delta_2$ satisfy $(x_0,y_0)\notin D$.
Define stopping time $\tau_D$ by
 \beqnn
\tau_D:=\inf\{t\ge0:(x_t,y_t)\in D\}.
 \eeqnn
Suppose that there exists a function $g\in C^2((0,\infty)\times(0,\infty))$ satisfying the following conditions:
\begin{itemize}
\item[{\normalfont(i)}]$0<\sup_{u,v>0}|g(u,v)|<\infty$ and $g(x_0,y_0)>0$;
\item[{\normalfont(ii)}] $\sup_{\varepsilon_1<u<\delta_1,\varepsilon_2<v<\delta_2}|\mathcal{L}g(u,v)|<\infty$;
\item[{\normalfont(iii)}]
$g(z_1,u)=g(u,z_2)= g(z_3,u)=g(u,z_4)=0$ for all $u>0$, $z_1\ge \delta_1$, $z_2\ge \delta_2$, $z_3\le \varepsilon_1$ and $z_4\le \varepsilon_2$;
\item[{\normalfont(iv)}]
There is a constant $d>0$ so that
$\mathcal{L}g(u,v)\ge d g(u,v)$ for all $u\in(\varepsilon_1,\delta_1),v\in(\varepsilon_2,\delta_2)$ and $(u,v)\notin D$.
\end{itemize}
Then
$\mbf{P}\{\tau_D<\infty\}\ge g(x_0,y_0)/\sup_{u,v>0}|g(u,v)|$.
\eproposition
\proof
Define the stopping time $\tau$ by
$\tau:=\tau_D\wedge\tau_{\varepsilon_1}^-(x)\wedge\tau_{\varepsilon_2}^-(y)
\wedge\tau_{\delta_1}^+(x)\wedge\tau_{\delta_2}^+(y)$.
By \eqref{3.1},
 \beqlb\label{3.22}
M^g_{t\wedge\tau }:=g(x_{t\wedge\tau },
y_{t\wedge\tau })
-g(x_0,y_0)
-\int_0^{t\wedge\tau }\mathcal{L}g(x_s,y_s)\dd s
 \eeqlb
is a local martingale.
Let $c_1:=\sup_{u,v>0}|g(u,v)|$ and $c_2:=\sup_{\varepsilon_1<u<\delta_1,\varepsilon_2<v<\delta_2}|\mathcal{L}g(u,v)|$.
Then $c_1,c_2<\infty$ by conditions (i) and (ii).
Moreover, $|\mathcal{L}g(x_s,y_s)|\le c_2$ for all $s<\tau$ and then
 \beqnn
|M^g_{t\wedge\tau }|\le |g(x_{t\wedge\tau },
y_{t\wedge\tau })|
+|g(x_0,y_0)|
+\int_0^{t\wedge\tau }|\mathcal{L}g(x_{s-},y_{s-})|\dd s
\le 2c_1+tc_2.
 \eeqnn
 Thus, from \cite[p.38]{P05} it follows that $t\mapsto M^g_{t\wedge\tau }$
is a martingale.
Taking the expectation on both sides of \eqref{3.22} we get
 \beqnn
\mbf{E}\big[g(x_{t\wedge\tau },
y_{t\wedge\tau })\big]
=g(x_0,y_0)+\int_0^t\mbf{E}\big[\mathcal{L}g(x_s,y_s)
1_{\{s\le\tau \}}\big]\dd s.
 \eeqnn
Taking integration by parts, we have
 \beqnn
 \ar\ar
\e^{-d  t}\mbf{E}\big[g(x_{t\wedge\tau },y_{t\wedge\tau })\big] \cr
 \ar=\ar
g(x_0,y_0)
-d \int_0^t\e^{-d s}\mbf{E}\big[g(x_{s\wedge\tau },
y_{s\wedge\tau })\big]\dd s
+\int_0^t\e^{-d  s}\dd\Big(\mbf{E}\big[g(x_{s\wedge\tau },
y_{s\wedge\tau })\big]\Big) \cr
 \ar=\ar
g(x_0,y_0)
-d \int_0^t\e^{-d  s}\mbf{E}\big[g(x_{s\wedge\tau },
y_{s\wedge\tau })\big]\dd s
+\int_0^t\e^{-d  s} \mbf{E}\big[\mathcal{L}g(x_s,y_s)
1_{\{s\le\tau \}}\big]
\dd s.
 \eeqnn
Now we conclude by  condition (iv) that
 \beqnn
\e^{-d  t}\mbf{E}\big[g(x_{t\wedge\tau },y_{t\wedge\tau })\big]
\ge
g(x_0,y_0)
-d \int_0^t\e^{-d s}\mbf{E}\big[g(x_{\tau },y_{\tau })
1_{\{s> \tau \}}\big]\dd s.
 \eeqnn
Leaving $t\to\infty$ and using condition (i),
 \beqnn
g(x_0,y_0)
 \ar\le\ar
d \int_0^\infty\e^{-d s}\mbf{E}\big[g(x_{\tau },y_{\tau })
1_{\{s\ge \tau \}}\big]\dd s
=
\mbf{E}\big[g(x_{\tau },y_{\tau })
\e^{-d  \tau } 1_{\{\tau<\infty\}}\big] \cr
 \ar\le\ar
\mbf{E}\Big[|g(x_{\tau },y_{\tau })|
\e^{-d  \tau }
\big[1_{\{\tau_{\varepsilon_1}^-(x)\wedge\tau_{\varepsilon_2}^-(y)
\le\tau_D\wedge\tau_{\delta_1}^+(x)\wedge\tau_{\delta_2}^+(y),
\tau_{\varepsilon_1}^-(x)\wedge\tau_{\varepsilon_2}^-(y)<\infty\}} \cr
 \ar\ar
\qquad\qquad\qquad\quad
+1_{\{\tau_{\delta_1}^+(x)\wedge\tau_{\delta_2}^+(y)
\le\tau_D\wedge\tau_{\varepsilon_1}^-(x)\wedge\tau_{\varepsilon_2}^-(y),
\tau_{\delta_1}^+(x)\wedge\tau_{\delta_2}^+(y)<\infty\}} \cr
 \ar\ar
\qquad\qquad\qquad\quad
+1_{\{\tau_D\le\tau_{\varepsilon_1}^-(x)\wedge\tau_{\varepsilon_2}^-(y)
\wedge\tau_{\delta_1}^+(x)\wedge\tau_{\delta_2}^+(y),
\tau_D<\infty \}}\big]\Big] \cr
 \ar\le\ar
\mbf{E}\Big[|g(x_{\tau },y_{\tau })|
\e^{-d  \tau }
1_{\{\tau_{\varepsilon_1}^-(x)\wedge\tau_{\varepsilon_2}^-(y)
\le\tau_D\wedge\tau_{\delta_1}^+(x)\wedge\tau_{\delta_2}^+(y),
\tau_{\varepsilon_1}^-(x)\wedge\tau_{\varepsilon_2}^-(y)<\infty\}}\Big] \cr
 \ar\ar
+\mbf{E}\Big[|g(x_{\tau },y_{\tau })|
\e^{-d  \tau }1_{\{\tau_{\delta_1}^+(x)\wedge\tau_{\delta_2}^+(y)
\le\tau_D\wedge\tau_{\varepsilon_1}^-(x)\wedge\tau_{\varepsilon_2}^-(y),
\tau_{\delta_1}^+(x)\wedge\tau_{\delta_2}^+(y)<\infty\}}\Big] \cr
 \ar\ar
+\sup_{u,v>0}|g(u,v)|\mbf{P}\{\tau_D<\tau_{\varepsilon_1}^-(x)
\wedge\tau_{\varepsilon_2}^-(y)
\wedge\tau_{\delta_1}^+(x)\wedge\tau_{\delta_2}^+(y),\tau_D<\infty\} \cr
 \ar=:\ar
I_1+I_2+I_3,
 \eeqnn
Since $t\mapsto x_t$ and $t\mapsto y_t$ are right continuous,
 under condition (iii), $I_1=I_2=0$.
Observe that $I_3\le
\sup_{u,v>0}|g(u,v)|\mbf{P}\{\tau_D<\infty\}$.
The proof is concluded.
\qed

The following result concerns a criterion for
the probability of extinction to be strictly less than one.

\bproposition\label{Pro2.1}
Suppose that
\begin{itemize}
\item[{\normalfont(i)}]
there are a constant $\beta>0$, a nonnegative function $g\in C^2((0,\infty)\times(0,\infty))$ and a strictly positive and
non-increasing function $\bar{g}$ on $(0,\infty)$ so that
$g(u,v)\ge \bar{g}(vu^{-\beta})$ for all $u,v>0$.
\item[{\normalfont(ii)}]
there are constants $0<\varepsilon_0<1$ and $u_*>0$
such that 
 \beqlb\label{10.9}
\mathcal{L}g(u,v)\le 0,\qquad  0<u,v\le\varepsilon_0,~~ vu^{-\beta}>u_*
 \eeqlb
and
 \beqlb\label{10.9b}
\inf_{z_1\le u,v\le z_2}g(u,v)>0~~\mbox{and }~
\sup_{z_1\le u,v\le z_2}|\mathcal{L}g(u,v)|<\infty,\qquad \mbox{for all }
z_2>z_1>0;
 \eeqlb
 \item[{\normalfont(iii)}]
$\sup_{t\ge0} (x_t+y_t)<\infty$, $\tau_0^-(x)=\infty$ a.s.
and there are constants $C,\delta>0$ independent of $0<\varepsilon\le\varepsilon_0$
such that
 \beqlb\label{10.3a}
\mbf{P}\Big\{\sup_{t\ge0}(x_t+y_t)\ge\varepsilon\Big\}
\le C\varepsilon^{\delta}.
 \eeqlb
\end{itemize}
In addition, there exist constants $u_0> u_*$ and $\tilde{\varepsilon}>0$ such that $g(x_0,y_0)<(1-\tilde{\varepsilon})\bar{g}(u_*)$
for all $(x_0,y_0)=(\varepsilon^{1+\beta^{-1}}, u_0\varepsilon^{\beta+1})$ with $0<\varepsilon<\varepsilon_0$. Then there is a constant $0<\varepsilon_1\le\varepsilon_0$ so that
 $\mbf{P}\{\tau_0^-(y)<\infty\}<1$ for the above $(x_0, y_0)$ with $0<\varepsilon<\varepsilon_1$.
\eproposition

Proposition \ref{Pro2.1} provides a criterion
for non-extinction with positive probability in which functions $g$ and $\bar{g}$ are chosen to allow test functions involving the ratio of $y$ and $x^\beta$.
It will be applied to prove Theorem \ref{t0.02}
for  functions $g(u,v):=u^{\beta\delta}v^{-\delta}+v^{\rho}$
and $\bar{g}(u):=u^{-\delta}$, which satisfy all the assumptions of Proposition \ref{Pro2.1}.
The key assumption of \eqref{10.9} is verified in Lemma \ref{t2.5},
and the other assumptions are verified in Lemma \ref{t10.1};
see the details in Subsection \ref{subsection3.3}.
The key to proving Proposition \ref{Pro2.1} is considering
a ratio $U$ of processes $y$ and $x^\beta$ bounded away from small constant $u_*>0$.

{\it Proof of Proposition \ref{Pro2.1}.}
Define the process $U:=(U_t)_{t\ge0}$ by $U_t=y_tx_t^{-\beta}$.
Since $x_0=\varepsilon^{1+\beta^{-1}}$ and $y_0=u_0\varepsilon^{\beta+1}$,
then $U_0=x_0^{-\beta}y_0=u_0$.
Under the assumption of $g(x_0,y_0)<(1-\tilde{\varepsilon})\bar{g}(u_*)$,
there is a constant $0<\varepsilon_1<\varepsilon_0$ so that for all $0<\varepsilon\le \varepsilon_1$ we have
 \beqlb\label{10.4}
1-g(x_0,y_0)/\bar{g}(u_*)- C\varepsilon^{\delta}>0,
 \eeqlb
where the constant $C$ given in \eqref{10.3a} does not depend on   $\varepsilon>0$.
Using \eqref{3.1},
 \beqnn
M_{t\wedge\gamma_{m,n}}^g
:=
g(x_{t\wedge\gamma_{m,n}},y_{t\wedge\gamma_{m,n}})-g(x_0,y_0)
-\int_0^{t\wedge\gamma_{m,n}}\mathcal{L}g(x_s,y_s)\dd s
 \eeqnn
is a local martingale, where $\gamma_{m,n}:=\tau_{m,n}\wedge
\tau_{u_*}^-(U)$.
Let
 \beqnn
H_s:=g(x_s,y_s)^{-1}\mathcal{L}g(x_s,y_s).
 \eeqnn
Taking integration by parts,
 \beqlb\label{10.5}
 \ar\ar
g(x_{t\wedge\gamma_{m,n}},y_{t\wedge\gamma_{m,n}})
\e^{-\int_0^{t\wedge\gamma_{m,n}} H_s\dd s} \cr
 \ar=\ar
g(x_0,y_0)
+\int_0^tg(x_{s\wedge\gamma_{m,n}},y_{s\wedge\gamma_{m,n}})\dd \big(\e^{-\int_0^{s\wedge\gamma_{m,n}} H_v\dd v}\big)
+\int_0^t \e^{-\int_0^{s\wedge\gamma_{m,n}} H_v\dd v} \dd g(x_{s\wedge\gamma_{m,n}},y_{s\wedge\gamma_{m,n}}) \cr
 \ar=\ar
g(x_0,y_0)
-\int_0^tg(x_{s\wedge\gamma_{m,n}},y_{s\wedge\gamma_{m,n}})
\e^{-\int_0^{s\wedge\gamma_{m,n}} H_v\dd v}
H_s1_{\{s\le\gamma_{m,n}\}}\dd s \cr
 \ar\ar
+\int_0^t \e^{-\int_0^{s\wedge\gamma_{m,n}} H_v\dd v}
\mathcal{L}g(x_s,y_s)1_{\{s\le\gamma_{m,n}\}} \dd s
+W_{m,n}(t)  \cr
 \ar=\ar
g(x_0,y_0)+W_{m,n}(t),
 \eeqlb
where $W_{m,n}(t):=\int_0^t \e^{-\int_0^{s\wedge\gamma_{m,n}} H_v\dd v}
\dd M^g_{s\wedge\gamma_{m,n}}$.
By \eqref{10.9b}, $\sup_{v\ge0}|H_v1_{\{v<\gamma_{m,n}\}}|\le C_{m,n}$
for some constant $C_{m,n}>0$. It follows from \cite[p.128]{P05} that
$t\mapsto W_{m,n}(t)$ is a local martingale.
Then there are stopping times $\gamma_k$ so that $\gamma_k\to\infty$ as $k\to\infty$
and $t\mapsto W_{m,n}(t\wedge\gamma_k)$
is a martingale for each $k\ge1$.
Let $\gamma_{m,n,k}:=\gamma_{m,n}\wedge\gamma_k$.
It then follows from \eqref{10.5} that
 \beqnn
g(x_0,y_0)=\mbf{E}\Big[g(x_{t\wedge\gamma_{m,n,k}},y_{t\wedge\gamma_{m,n,k}})
\e^{-\int_0^{t\wedge\gamma_{m,n,k}} H_s\dd s}\Big].
 \eeqnn
By the condition (iii), we have $\tau_0^-(x)=\tau_\infty^+(x)=\tau_\infty^+(y)=\infty$
a.s. and then $x_t>0$ for all $t\ge0$ a.s.
Then
 \beqlb\label{10.2}
\tau_0^-(y)=\tau_0^-(U)\ge \tau_{u_*}^-(U)
 \eeqlb
almost surely.  Moreover,
 \beqnn
\lim_{m,n,k\to\infty}\gamma_{m,n,k}
=\tau_0^-(x)\wedge\tau_0^-(y)\wedge\tau_\infty^+(x)\wedge\tau_\infty^+(y)
\wedge\tau_{u_*}^-(U)
=\tau_{u_*}^-(U)
 \eeqnn
almost surely.
It thus follows from Fatou's lemma that
 \beqnn
g(x_0,y_0)
 \ar=\ar
\liminf_{t,m,n,k\to\infty}\mbf{E}\Big[g(x_{t\wedge\gamma_{m,n,k}},y_{t\wedge\gamma_{m,n,k}})
\e^{-\int_0^{t\wedge\gamma_{m,n,k}}H_s\dd s}\Big] \cr
 \ar\ge\ar
\mbf{E}\Big[\liminf_{t,m,n,k\to\infty}g(x_{t\wedge\gamma_{m,n,k}},y_{t\wedge\gamma_{m,n,k}})
\e^{-\int_0^{t\wedge\gamma_{m,n,k}}H_s\dd s}\Big] \cr
 \ar\ge\ar
\mbf{E}\Big[g(x_{\tau_{u_*}^-(U)},y_{\tau_{u_*}^-(U)})
1_{\{\tau_{u_*}^-(U)<\infty\}}
\e^{-\int_0^{\tau_{u_*}^-(U)}H_s\dd s}
\Big].
 \eeqnn
Since $U_{\tau_{u_*}^-(U)}\le u_*$
for $\tau_{u_*}^-(U)<\infty$,
then $\bar{g}(U_{\tau_{u_*}^-(U)})\ge \bar{g}(u_*)$
since $\bar{g}$ is non-increasing when $\tau_{u_*}^-(U)<\infty$.
It then under condition (i),
 \beqnn
g(x_0,y_0)
\ge
\bar{g}(u_*)\mbf{E}\Big[
1_{\{\tau_{u_*}^-(U)<\infty\}}
\e^{-\int_0^{\tau_{u_*}^-(U)}H_s\dd s}
\Big].
 \eeqnn
For $0<\varepsilon<\varepsilon_1$ and for $s< \tau_{u_*}^-(U)$ and $\sup_{t\ge0}(x_t+y_t)\le \varepsilon$, we have $H_s\le0$ by \eqref{10.9}.
Thus
 \beqnn
g(x_0,y_0)/\bar{g}(u_*)
 \ar\ge\ar
\mbf{E}\Big[
\e^{-\int_0^{\tau_{u_*}^-(U)}H_v\dd v}
1_{\{\sup_{s\ge0}(x_s+y_s)\le \varepsilon,\tau_{u_*}^-(U)<\infty\}}\Big] \cr
 \ar\ge\ar
\mbf{P}\Big\{\sup_{s\ge0}(x_s+y_s)\le \varepsilon,\tau_{u_*}^-(U)<\infty\Big\}
 =:
M_\varepsilon.
 \eeqnn
Applying \eqref{10.3a},
 \beqnn
\mbf{P}\big\{\tau_{u_*}^-(U)<\infty\big\}
 \ar=\ar
M_\varepsilon
+\mbf{P}\Big\{\sup_{s\ge0}(x_s+y_s)> \varepsilon,\tau_{u_*}^-(U)<\infty\Big\} \cr
 \ar\le\ar
M_\varepsilon
+\mbf{P}\Big\{\sup_{s\ge0}(x_s+y_s)> \varepsilon\Big\}
\le
g(x_0,y_0)/\bar{g}(u_*)+C\varepsilon^{\delta}.
 \eeqnn
It then follows from \eqref{10.4} that
 \beqnn
\mbf{P}\big\{\tau_{u_*}^-(U)=\infty\big\}
\ge
1-g(x_0,y_0)/\bar{g}(u_*)- C\varepsilon^{\delta}>0.
 \eeqnn
Then $\mbf{P}\{\tau_0^-(y)=\infty\}>0$ follows by \eqref{10.2}.
\qed

\section{Proofs of Theorems \ref{t0.01}--\ref{t0.03}}\label{proof of Theorems}

\setcounter{equation}{0}

In this section, we use the propositions and corollary in Section \ref{Criteria}
to complete the proof of Theorems \ref{t0.01}--\ref{t0.03}.
We always assume that $((X_t,Y_t))_{t\ge0}$ is a solution to \eqref{a1.1}.
For any $x,y,z\ge0$ and $g\in C^2((0,\infty)\times(0,\infty))$ define
 \beqlb\label{3.3}
K_z^1g(x,y):=
g(x+z,y)-g(x,y)-g'_x(x,y)z
 \eeqlb
and
 \beqlb\label{3.4}
K_z^2g(x,y):=
g(x,y+z)-g(x,y)-g'_y(x,y)z.
 \eeqlb
Let
 \beqlb\label{10.1}
\mathcal{L}g(x,y)
 \ar:=\ar
-a_1x^{p_1+1}g_x'(x,y)
+a_2x^{p_2+2}g_{xx}''(x,y)
+a_3x^{p_3+\alpha_1}\int_0^\infty K_z^1g(x,y)\mu_1(\dd z) \cr
 \ar\ar
-b_1y^{q_1+1}g_y'(x,y)
+b_2y^{q_2+2}g_{yy}''(x,y)
+b_3y^{q_3+\alpha_2}\int_0^\infty K_z^2g(x,y)\mu_2(\dd z) \cr
 \ar\ar
-\eta_1x^{\theta_1}y^{\kappa_1}g_x'(x,y)
-\eta_2y^{\theta_2}x^{\kappa_2}g_y'(x,y).
 \eeqlb
For any $g\in C^2((0,\infty))$ and $x,z>0$ define
 \beqlb\label{3.2}
K_zg(x):=g(x+z)-g(x)-g'(x)z.
 \eeqlb
Substituting the variable $z$  by $zx$ we get
 \beqlb\label{3.5}
\int_0^\infty g(z)\mu_i(\dd z)=
x^{-\alpha_i}\int_0^\infty g(zx)\mu_i(\dd z),\qquad i=1,2,
 \eeqlb
where we use the fact $\mu_i(\dd z)=\frac{\alpha_i(\alpha_i-1)}{\Gamma(\alpha_i)\Gamma(2-\alpha_i)}z^{-1-\alpha_i}\dd z$.
By Taylor's formula, for any bounded function $g$ with continuous second derivative,
 \beqlb\label{4.1a}
g(x+z)-g(x)=z\int_0^1g'(x+zv)\dd v,
 \eeqlb
and
 \beqlb\label{4.1}
K_zg(x)=z^2\int_0^1g''(x+zv)(1-v)\dd v.
 \eeqlb

\subsection{Preliminaries}

In this subsection, we present some preliminary results that are needed
for the proof of Theorems \ref{t0.01}--\ref{t0.03}.

\blemma\label{t10.1}
(i) $\sup_{t\ge0} (X_t+Y_t)<\infty$ and $\tau_\infty^+(X)=\tau_\infty^+(Y)=\infty$ a.s.

(ii) Given $X_0,Y_0>0$, $\delta\in(0,1/2)$ and $c_0>0$,
there exists a constant $C>0$ independent of $X_0,Y_0,c_0$ satisfying
 \beqnn
\mbf{P}\Big\{\sup_{t\ge0}X_t\ge c_0\Big\}\le Cc_0^{-\delta}X_0^{\delta},\quad
\mbf{P}\Big\{\sup_{t\ge0}Y_t\ge c_0\Big\}\le
Cc_0^{-\delta}Y_0^{\delta}.
 \eeqnn
\elemma
\proof
(i) By \cite[Theorem 2.3(i) and Proposition 2.6]{LYZh},
$X_t\to0$ and $Y_t\to0$ as $t\to\infty$  a.s. Then $\sup_{t\ge0} (X_t+Y_t)<\infty$ a.s.

(ii)
We first estimate $\mbf{E}[X_t^\delta+Y_t^\delta]<\infty$ using the
same arguments as in
the proof of \cite[Lemma 3.3]{RXYZ19} and then use the Markov inequality to get the assertion.
\qed

In the following proposition, we present a comparison theorem for
the SDE system \eqref{a1.1}, which is used in the proof of Theorem \ref{t0.02}.
The proof for the following proposition is a modification
of the classical Yamada-Watanabe argument for the pathwise uniqueness of SDE
and is stated in Section \ref{Appendix}.

\bproposition\label{t2.1} (Comparison theorem)
If $(\tilde{X}_t,\tilde{Y}_t)_{t\ge0}$ is another solution to \eqref{a1.1}
satisfying $\tilde{X}_0\ge X_0$ and $\tilde{Y}_0\le Y_0$,
then we have $\mbf{P}\{\tilde{X}_t\ge X_t$ and $\tilde{Y}_t\le Y_t$ for all
$0\le t<\tau\}=1$, where
\[\tau:=\tau_0^-(X)\wedge\tau_0^-(Y)\wedge
\tau_0^-(\tilde{X})\wedge\tau_0^-(\tilde{Y}).\]
\eproposition

In the following, we prove an irreducibility result
for SDE system \eqref{a1.1}.

\bproposition\label{t2.10}
Given $x_2>x_1>0$, $y_2>y_1>0$ and $D=[x_1, x_2]\times[y_1,y_2]$
for $(X_0,Y_0)\notin D$ and $\tau_D:=\inf\{t\ge0:(X_t,Y_t)\in D\}$ we have
 $\mbf{P}\{\tau_D<\infty\}>0$.
\eproposition

The proof of Proposition \ref{t2.10} is an application of Proposition \ref{pro2.2a}
where the key is to choose the function $g$ of the form $h_{\lambda,\lambda_1}$ in Lemma 3.4.
Recall the operator $K_z$ defined in \eqref{3.2}.
To confirm condition (iv) of Proposition \ref{pro2.2a},
we need to estimate $h_{\lambda,\lambda_1}',h_{\lambda,\lambda_1}''$ and the integral $\int_0^\infty  K_zh_{\lambda,\lambda_1}(x)z^{-1-\alpha}\dd z$
in the following lemma.

\blemma\label{t2.4}
Given $0<x_1<x_2<x_3$, $\alpha\in(1,2)$
and
\[h_{\lambda,\lambda_1}(x)=\e^{-\lambda/(x-x_1)-\lambda\lambda_1/(x_3-x)}1_{\{x_1<x<x_3\}},\qquad x,\lambda,\lambda_1>0,\]
we have
\begin{itemize}
\item[{\normalfont(i)}]
$h_{\lambda,\lambda_1}'(x)\le\lambda h_{\lambda,\lambda_1}(x)(x-x_1)^{-2}$ for all $x_1<x<x_3$ and $\lambda,\lambda_1>0$.

\item[{\normalfont(ii)}]
there are constants $\lambda_0,\lambda_1,c_0>0$ so that for all $\lambda>\lambda_0$ and all $x_1<x<x_2$, we have
 \beqlb\label{2.5a}
h_{\lambda,\lambda_1}''(x)h(x)^{-1}\ge\lambda^2 c_0(x-x_1)^{-4}
 \eeqlb
and
 \beqlb\label{2.5b}
h_{\lambda,\lambda_1}(x)^{-1} \int_0^\infty K_zh_{\lambda,\lambda_1}(x)z^{-1-\alpha}\dd z
 \ge
\lambda^\alpha c_0(x-x_1)^{-2-\alpha}.
 \eeqlb
\item[{\normalfont(iii)}]
for each $\lambda_1>0$ there are constants $\lambda_0,\tilde{c}_0,\tilde{c}_1>0$ so that for all $\lambda\ge\lambda_0$ and all $x_1<x<x_3$, we have
 \beqlb\label{2.5c}
h_{\lambda,\lambda_1}''(x)h_{\lambda,\lambda_1}(x)^{-1}\ge\lambda^2 \tilde{c}_0\big[(x-x_1)^{-4}
-\tilde{c}_1\big]
 \eeqlb
and
 \beqlb\label{2.5d}
h_{\lambda,\lambda_1}(x)^{-1} \int_0^\infty  K_zh_{\lambda,\lambda_1}(x)z^{-1-\alpha}\dd z
\ge
\lambda^{\alpha} \tilde{c}_0\big[(x-x_1)^{-2-\alpha}
-\tilde{c}_1\big].
 \eeqlb
\end{itemize}
\elemma
\proof
For simplicity, we assume that $x_1=1,x_2=2,x_3=3$.
Observe that for all $1<x<3$,
 \beqlb\label{2.3}
h_{\lambda,\lambda_1}'(x)=\lambda h_{\lambda,\lambda_1}(x)\big[(x-1)^{-2}-\lambda_1(x-3)^{-2}\big].
 \eeqlb
Then assertion (i) follows.
For all $1<x<3$,
 \beqlb\label{2.3b}
\lambda^{-1}h_{\lambda,\lambda_1}''(x)h_{\lambda,\lambda_1}(x)^{-1}
 \ar=\ar
\lambda(x-1)^{-4}+\lambda\lambda_1^2(x-3)^{-4} \cr
 \ar\ar
-2\lambda\lambda_1(x-1)^{-2}(x-3)^{-2}
-2(x-1)^{-3}+2\lambda_1(x-3)^{-3}.
 \eeqlb
Observe that
 \beqlb\label{2.3a}
2^{-1}\lambda(x-1)^{-4}>2(x-1)^{-3},~~
2^{-1}\lambda\lambda_1^2(x-3)^{-4}>-2\lambda_1(x-3)^{-3},
\qquad 1<x<3
 \eeqlb
if $\lambda,\lambda\lambda_1>8$.
For all $1<x<5/2$, we have $[(x-3)/(x-1)]^2>1/9$
and then for $\lambda_1<72^{-1}$, we have $\lambda_1<(4+\sqrt{14})^{-1}9^{-1}<(4+\sqrt{14})^{-1}[(x-3)/(x-1)]^2$,
which implies
 \beqnn
2\lambda_1^2
-8\lambda_1[(x-3)/(x-1)]^{2}+[(x-3)/(x-1)]^{4}\ge0,\qquad 1<x<5/2.
 \eeqnn
It follows from \eqref{2.3b} and \eqref{2.3a} that
for $\lambda,\lambda\lambda_1>8$ and $\lambda_1<72^{-1}$,
 \beqlb\label{2.1}
 \ar\ar
h_{\lambda,\lambda_1}''(x)h_{\lambda,\lambda_1}(x)^{-1} \cr
 \ar=\ar
4^{-1}\lambda^2(x-1)^{-4}
+2^{-1}\lambda^2(x-1)^{-4}-2\lambda(x-1)^{-3}
+2^{-1}\lambda^2\lambda_1^2(x-3)^{-4}
+2\lambda\lambda_1(x-3)^{-3} \cr
 \ar\ar
+4^{-1}\lambda^2(x-3)^{-4}\big[[(x-3)/(x-1)]^{4}+2\lambda_1^2
-8\lambda_1[(x-3)/(x-1)]^{2}\big] \cr
 \ar\ge\ar
4^{-1}\lambda^2(x-1)^{-4}
,\qquad 1<x<5/2,
 \eeqlb
which gives \eqref{2.5a}.
For $0<z<(\lambda\lambda_1)^{-1}\wedge2^{-1}$ and $1<x<2$, we have
$\lambda\lambda_1 z(3-x-z)^{-1}(3-x)^{-1}\le2$ and then
$h_{\lambda,\lambda_1}(x+z)\ge h_{\lambda,\lambda_1}(x)\e^{-2}$. Thus, by \eqref{4.1} and \eqref{2.1}, for all
for $\lambda,\lambda\lambda_1>8$ and $\lambda_1<72^{-1}$,
 \beqnn
 \ar\ar
\int_{(\lambda\lambda_1)^{-1}}^{1/2} K_zh_{\lambda,\lambda_1}(x)z^{-1-\alpha}\dd z
=
\int_{(\lambda\lambda_1)^{-1}}^{1/2}h_{\lambda,\lambda_1}''(x+zv)z^{1-\alpha}\dd z\int_0^1(1-v)\dd v \ge0,~~ 1<x<2
 \eeqnn
and
 \beqnn
\int_0^{(\lambda\lambda_1)^{-1}} K_zh_{\lambda,\lambda_1}(x) z^{-1-\alpha}\dd z
 \ar=\ar
\int_0^{(\lambda\lambda_1)^{-1}}h_{\lambda,\lambda_1}''(x+zv)z^{1-\alpha}\dd z\int_0^1(1-v)\dd v \cr
 \ar\ge\ar
h_{\lambda,\lambda_1}(x)\e^{-2}8^{-1}\lambda^2\int_0^{(\lambda\lambda_1)^{-1}}(x+z-1)^{-4}
z^{1-\alpha}\dd z.
 \eeqnn
For $\lambda\lambda_1\ge1$ and $1<x<2$, with the variable $z$ replaced by
$(x-1)u$ we get
 \beqnn
 \ar\ar
\int_0^{(\lambda\lambda_1)^{-1}}(x+z-1)^{-4}
z^{1-\alpha}\dd z
 =
(x-1)^{-2-\alpha}
\int_0^{[\lambda\lambda_1(x-1)]^{-1}}(1+u)^{-4}
u^{1-\alpha}\dd u \cr
 \ar\ge\ar
(x-1)^{-2-\alpha}2^{-4}
\int_0^{(\lambda\lambda_1)^{-1}}
u^{1-\alpha}\dd u
=
(x-1)^{-2-\alpha}2^{-4}(2-\alpha)^{-1}
(\lambda\lambda_1)^{\alpha-2}.
 \eeqnn
It follows that for $\lambda,\lambda\lambda_1>8$ and $\lambda_1<72^{-1}$,
\beqlb\label{2.4}
\int_0^{1/2} K_zh_{\lambda,\lambda_1}(x) z^{-1-\alpha}\dd z
\ge
2c_0\lambda^{\alpha} (x-1)^{-2-\alpha}h_{\lambda,\lambda_1}(x),
\quad 1<x<2.
 \eeqlb
where $c_0:=\e^{-2}2^{-8}(2-\alpha)^{-1}\lambda_1^{\alpha-2}$.
 Using \eqref{2.3},
 \beqnn
h_{\lambda,\lambda_1}(x)^{-1}\int_{1/2}^\infty K_zh_{\lambda,\lambda_1}(x)z^{-1-\alpha}\dd z
 \ar\ge\ar
-\int_{1/2}^\infty z^{-1-\alpha}\dd z
-(x-1)^{-2}\lambda \int_{1/2}^\infty z^{-\alpha}\dd z \cr
 \ar\ge\ar
-c_0\lambda^{\alpha} (x-1)^{-2-\alpha},\qquad 1<x<2
 \eeqnn
for all large enough $\lambda>0$.
Therefore, \eqref{2.5b} follows from \eqref{2.4}.

In the following, we prove  assertion (iii) and take $\lambda_1=1$.
By \eqref{2.3b} and \eqref{2.3a}, there are constants $\delta\in(0,1/2)$ and $\tilde{c}_0>0$ such that for all large enough $\lambda>1$,
 \beqlb\label{2.7}
\lambda^{-2}h_{\lambda,1}''(x)h_{\lambda,1}(x)^{-1}
\ge
\tilde{c}_0(x-1)^{-4},
\quad x\in(1,1+2\delta)\cup(3-2\delta,3).
 \eeqlb
By \eqref{2.3b}, there is a constant $\tilde{c}_1>0$ such that
for all $\lambda>1$
 \beqnn
\lambda^{-2}h_{\lambda,1}''(x)h_{\lambda,1}(x)^{-1}
 \ar\ge\ar
-2(x-1)^{-2}(x-3)^{-2}
-2\lambda^{-1}(3-x)^{-3} \cr
 \ar\ge\ar
-2(x-1)^{-2}(x-3)^{-2}
-2(3-x)^{-3} \cr
 \ar=\ar
\tilde{c}_0(x-1)^{-4}
-[\tilde{c}_0(x-1)^{-4}
+2(x-1)^{-2}(x-3)^{-2}
+2(3-x)^{-3}] \cr
 \ar\ge\ar
\tilde{c}_0(x-1)^{-4}-\tilde{c}_1, \qquad x\in[1+2\delta,3-2\delta].
 \eeqnn
Thus, \eqref{2.5c} follows from \eqref{2.7}.
For $x\in(1,1+\delta)\cup(3-2\delta,3)$, by \eqref{2.7} and the same arguments as
in \eqref{2.4}, we have
 \beqnn
\int_0^{\delta}K_zh_{\lambda,1}(x)z^{-1-\alpha}\dd z
\ge
\bar{c}_0\lambda^{\alpha}h_{\lambda,1}(x)(x-1)^{-2-\alpha}
 \eeqnn
for some constant $\bar{c}_0>0$.
For $1+\delta\le x\le 3-\delta$ and $0<z\le \delta$ let
 \beqnn
\tilde{h}(x,z):=1/(x+z-1)-1/(3-x-z)-[1/(x-1)-1/(3-x)].
 \eeqnn
Then there is a constant $\hat{c}_0>0$ so that
$|\tilde{h}''_{zz}(x,u)|\le \hat{c}_0$ for $1+\delta\le x\le 3-\delta$ and $0<u\le \delta$.
Since $\e^u-1\ge u$ for all $u\in\mbb{R}$,
then $h_{\lambda,1}(x)^{-1}K_zh_{\lambda,1}(x)\ge-\lambda[\tilde{h}(x,z)
+\tilde{h}'_z(x,0)]$
for $1+\delta\le x\le 3-\delta$ and $0<z\le \delta$.
Then by  \eqref{4.1},
 \beqnn
 \ar\ar
h_{\lambda,1}(x)^{-1}\int_0^{\delta}K_zh_{\lambda,1}(x)z^{-1-\alpha}\dd z
 \ge
-\lambda\int_0^{\delta}z^{1-\alpha}\dd z\int_0^1 (1-u)
\tilde{h}''_{zz}(x,zu)\dd u \cr
 \ar\ar\qquad\ge
-2^{-1}\hat{c}_0\lambda\int_0^{\delta}z^{1-\alpha}\dd z
=
\lambda^\alpha \bar{c}_0 (x-1)^{-2-\alpha}
-\lambda^\alpha \bar{c}_0 (x-1)^{-2-\alpha}
-2^{-1}\hat{c}_0\lambda\int_0^{\delta}z^{1-\alpha}\dd z \cr
 \ar\ar\qquad\ge
\lambda^\alpha \bar{c}_0 (x-1)^{-2-\alpha}
-\lambda^\alpha \bar{c}_0 \delta^{-2-\alpha}
-2^{-1}\hat{c}_0\lambda\int_0^{\delta}z^{1-\alpha}\dd z,\qquad 1+\delta\le x\le 3-2\delta.
 \eeqnn
Thus, for all $1<x<3$,
 \beqnn
h_{\lambda,1}(x)^{-1}\int_0^{\delta}K_zh_{\lambda,1}(x)z^{-1-\alpha}\dd z
\ge
\lambda^{\alpha}\bar{c}_0(x-1)^{-2-\alpha}
-\lambda^\alpha \bar{c}_0 \delta^{-2-\alpha}
-2^{-1}\hat{c}_0\lambda\int_0^{\delta}z^{1-\alpha}\dd z.
 \eeqnn
By \eqref{2.3}, for all $1<x<3$,
 \beqnn
h_{\lambda,1}(x)^{-1}\int_{\delta}^\infty K_zh_{\lambda,1}(x)z^{-1-\alpha}\dd z
\ge
- \int_{\delta}^\infty z^{-1-\alpha}\dd z
-(x-1)^{-2}\lambda\int_{\delta}^\infty z^{-\alpha}\dd z.
 \eeqnn
Thus, \eqref{2.5d} holds for all large enough $\lambda>0$.
\qed

{\it Proof of Proposition \ref{t2.10}.}
We first assume that $X_0\in(x_1,x_2)$, $Y_0\in (0,y_1)$ and $y_2\in (0,Y_0)$.
For $\lambda,\tilde{\lambda},\lambda_1,\tilde{\lambda}_1,x,y>0$ define
 \beqnn
g_{\lambda,\lambda_1,1}(x)=\e^{-\lambda/(x-x_1)-\lambda\lambda_1/(x_2-x)}1_{\{x_1<x<x_2\}},\quad
 g_{\tilde{\lambda},
 \tilde{\lambda}_1,2}(y)=\e^{-\tilde{\lambda}/(y-y_2)
 -\tilde{\lambda}\tilde{\lambda}_1/(y_3-y)}
 1_{\{y_2<y<y_3\}}.
 \eeqnn
Then $g_1,g_2\in C^2(0,\infty)$.
Let $g_{\lambda,\lambda_1,\tilde{\lambda},\tilde{\lambda}_1}(x,y)
=g_{\lambda,\lambda_1,1}(x)g_{\tilde{\lambda},\tilde{\lambda}_1,2}(y)$ for $x,y>0$ and we use Proposition \ref{pro2.2a} to establish the assertion with
$(\varepsilon_1,\delta_1,\varepsilon_2,\delta_2)=(x_1,x_2,y_2,y_3)$.
 Conditions (i)--(iii) in Proposition \ref{pro2.2a} are given by the definitions of $g_{\lambda,\lambda_1,1},g_{\tilde{\lambda},\tilde{\lambda}_1,2}$.
We verify condition (iv) in the following way.
By \eqref{10.1} and Lemma \ref{t2.4}, there are constants $\lambda_0,\tilde{\lambda}_0,\lambda_1,\tilde{\lambda}_1>0$ and
$c_0,c_1,\tilde{c}_0>0$ such that for all $\lambda>\lambda_0,\tilde{\lambda}>\tilde{\lambda}_0$
 \beqlb\label{2.10}
 \ar\ar
g_{\lambda,\lambda_1,\tilde{\lambda},\tilde{\lambda}_1}(x,y)^{-1}
\mathcal{L}g_{\lambda,\lambda_1,\tilde{\lambda},\tilde{\lambda}_1}(x,y) \cr
 \ar\ge\ar
-\lambda a_1(x-x_1)^{-2}x^{p_1+1}
-\lambda\eta_1(x-x_1)^{-2}x^{\theta_1}y^{\kappa_1}
+\lambda^2 c_0a_2\big[(x-x_1)^{-4}-c_1\big]x^{p_2+2} \cr
 \ar\ar
+\lambda^{\alpha_1} c_0a_3\big[(x-x_1)^{-2-\alpha_1}-c_1\big]x^{p_3+\alpha_1}
-\tilde{\lambda} b_1(y-y_2)^{-2}y^{q_1+1}
-\tilde{\lambda}\eta_2(y-y_2)^{-2}y^{\theta_2}x^{\kappa_2} \cr
 \ar\ar
+\tilde{\lambda}^2 \tilde{c}_0b_2 (y-y_2)^{-4} y^{q_2+2}
+\tilde{\lambda}^{\alpha_2} \tilde{c}_0b_3 (y-y_2)^{-2-\alpha_2} y^{q_3+\alpha_2} \cr
 \ar\ge\ar
I_{11}(\lambda,x)-I_{12}(\lambda)+I_2(\tilde{\lambda},y),
\qquad x_1<x<x_2,~y_2<y<y_1,
 \eeqlb
where $I_{12}(\lambda):=c_0c_1 [\lambda^2a_2x_2^{p_2+2}+\lambda^{\alpha_1}a_3x_2^{p_3+\alpha_1}]$,
 \beqnn
I_{11}(\lambda,x)
 \ar:=\ar
c_0\big[\lambda^2 a_2x_1^{p_2+2}(x-x_1)^{-4}+\lambda^{\alpha_1} a_3x_1^{p_3+\alpha_1}(x-x_1)^{-2-\alpha_1}\big] \cr
 \ar\ar
-\lambda a_1(x-x_1)^{-2}x_2^{p_1+1}
-\lambda\eta_1(x-x_1)^{-2}x_2^{\theta_1}y_1^{\kappa_1}
 \eeqnn
and
 \beqnn
I_2(\tilde{\lambda},y)
 \ar:=\ar
\tilde{c}_0\big[\tilde{\lambda}^2 b_2 (y-y_2)^{-4} y_2^{q_2+2}
+\tilde{\lambda}^{\alpha_2} b_3 (y-y_2)^{-2-\alpha_2} y_2^{q_3+\alpha_2}\big] \cr
 \ar\ar
-\tilde{\lambda} (y-y_2)^{-2}\big[b_1y_1^{q_1+1}+\eta_2y_1^{\theta_2}x_2^{\kappa_2}\big].
 \eeqnn
Since $a_2+a_3>0$, there is a large enough constant $\lambda_2>\lambda_0$
such that for all $x_1<x<x_2$,
 \beqlb\label{2.11}
I_{11}(\lambda_2,x)
 \ar=\ar
(x-x_1)^{-2}\big[\lambda_2^2 c_0 a_2x_1^{p_2+2}(x-x_1)^{-2}+\lambda_2^{\alpha_1} c_0a_3x_1^{p_3+\alpha_1}(x-x_1)^{-\alpha_1}\big]  \cr
 \ar\ar\qquad\qquad\quad
-\lambda_2 a_1x_2^{p_1+1}
-\lambda_2\eta_1x_2^{\theta_1}y_1^{\kappa_1}\big] \cr
 \ar\ge\ar
(x-x_1)^{-2}\big[\lambda_2^2 c_0 a_2x_1^{p_2+2}(x_2-x_1)^{-2}+\lambda_2^{\alpha_1} c_0a_3x_1^{p_3+\alpha_1}(x_2-x_1)^{-\alpha_1}\big]  \cr
 \ar\ar\qquad\qquad\quad
-\lambda_2 a_1x_2^{p_1+1}
-\lambda_2\eta_1x_2^{\theta_1}y_1^{\kappa_1}\big]  \cr
 \ar\ge\ar
(x-x_1)^{-2}\lambda_2(a_2+a_3)
\ge
\lambda_2(x_2-x_1)^{-2}(a_2+a_3).
 \eeqlb
Since $b_2+b_3>0$, there is a large enough constant $\tilde{\lambda}_2>\tilde{\lambda}_0$
such that for all $y_2<y<y_1$,
 \beqlb\label{2.12}
I_2(\tilde{\lambda}_2,y)
 \ar=\ar
(y-y_2)^{-2}\big[\tilde{\lambda}_2^2 \tilde{c}_0 b_2 (y-y_2)^{-2} y_2^{q_2+2}
+\tilde{\lambda}_2^{\alpha_2} \tilde{c}_0 b_3 (y-y_2)^{-\alpha_2} y_2^{q_3+\alpha_2} \cr
 \ar\ar\qquad\qquad\quad
-\tilde{\lambda}_2 b_1y_1^{q_1+1}-\tilde{\lambda}_2 \eta_2y_1^{\theta_2}x_2^{\kappa_2}\big] \cr
 \ar\ge\ar
(y-y_2)^{-2}\big[\tilde{\lambda}_2^2 \tilde{c}_0 b_2 (y_1-y_2)^{-2} y_2^{q_2+2}
+\tilde{\lambda}_2^{\alpha_2} \tilde{c}_0 b_3 (y_1-y_2)^{-\alpha_2} y_2^{q_3+\alpha_2} \cr
 \ar\ar\qquad\qquad\quad
-\tilde{\lambda}_2 b_1y_1^{q_1+1}-\tilde{\lambda}_2 \eta_2y_1^{\theta_2}x_2^{\kappa_2}\big] \cr
 \ar\ge\ar
(y-y_2)^{-2} \tilde{\lambda}_2(b_2+b_3)
\ge \tilde{\lambda}_2(y_1-y_2)^{-2}(b_2+b_3)
 \eeqlb
and
 \beqlb\label{2.13}
2^{-1}\tilde{\lambda}_2(y_1-y_3)^{-2}(b_2+b_3)\ge I_{12}(\lambda_2).
 \eeqlb
Combining \eqref{2.10}--\eqref{2.13} one obtains
 \beqnn
\mathcal{L}g_{\lambda_2,\lambda_1,\tilde{\lambda}_2,\tilde{\lambda}_1}(x,y)
\ge
2^{-1}\tilde{\lambda}_2(y_1-y_3)^{-2}(b_2+b_3)
g_{\lambda_2,\lambda_1,\tilde{\lambda}_2,\tilde{\lambda}_1}(x,y),
\quad x_1<x<x_2,~y_2<y<y_1.
 \eeqnn
Therefore, condition (iv) in Proposition \ref{pro2.2a} is confirmed
with $d:=2^{-1}\tilde{\lambda}_2(y_1-y_3)^{-2}(b_2+b_3)$.
Now $\mbf{P}\{\tau_D<\infty\}>0$ by Proposition \ref{pro2.2a}.
We can also obtain the assertion when $Y_0\in(y_1,y_2)$ and $X_0\in (0,x_1)$.
Using similar arguments, one can obtain the conclusion
when $(X_0,Y_0)\in((x_1,x_2)\times(y_2,\infty))\cup((x_2,\infty)\times(y_1,y_2))$.
Then, by the strong Markov property, we end the proof for general $X_0,Y_0>0$.
\qed

By Proposition \ref{pro2.2a}
and the proof of Proposition \ref{t2.10},
there is a positive lower bound for $\mbf{P}(\tau_D<\infty)$.

\subsection{Proof of Theorem \ref{t0.01}}

In this subsection, we apply Proposition \ref{t3.1b} to prove Theorem \ref{t0.01}.
We first verify the assertion \eqref{3.0} in the following lemma,
which is proved by introducing a function involving $-\ln(x+y^\beta)$ and
estimating the first time for $X_t+Y_t^\beta$ to reach below $k^{-1}$
for some $\beta>1$.

\blemma\label{t3.2}
If $\theta_1\vee\theta_2\ge1$, then
$\mbf{P}\{\tau_0^-(X)=\tau_0^-(Y)<\infty\}=0$.
\elemma
\proof
We assume that $\theta_1\ge1$ and $\beta>[(1-\theta_2)/\kappa_2]\vee1$.
Let $n\ge2$ and the nonnegative function $g_n\in C^2((0,\infty)\times(0,\infty))$ satisfy $g_n(x,y)=\ln (n+n^{\beta}+1)-\ln(x+y^\beta)$ when $0<x,y<n$,
$g_n(x,y)=0$ when $x\vee y>3n^\beta$,
and
 \beqlb\label{2.19}
 \ar\ar
c_0:=\sup_{(x,y)\in([1,3n^\beta+1]\times(0,\infty))\cup
((0,\infty)\times[1,3n^\beta+1])}
[g_n(x,y)+|\partial_xg_n(x,y)|+|\partial_{xx}g_n(x,y)| \cr
 \ar\ar\qquad\qquad\qquad\qquad\qquad\qquad\qquad\qquad
+|\partial_yg_n(x,y)|
+|\partial_{yy}g_n(x,y)|]<\infty,
 \eeqlb
where $\partial_xg_n(x,y):=\partial_xg_n(x,y)/\partial x$,
and $\partial_{xx}g_n(x,y),\partial_yg_n(x,y),\partial_{yy}g_n(x,y)$
are defined similarly.
In the following, we estimate $\mathcal{L}g_n$.
By Lemma \ref{t2.3}, there are constants $0<\delta_1<1\wedge(\kappa_1/\beta)$
and
$(1-\theta_2)/\beta<\delta_2<\kappa_2\wedge1$ so that for all $0<x,y\le n$ we have
 \beqnn
(x+y^\beta)^{-1} x^{\theta_1}y^{\kappa_1}
\le
x^{\theta_1-(1-\delta_1)}y^{\kappa_1-\delta_1\beta}
\le
n^{\theta_1-(1-\delta_1)+\kappa_1-\delta_1\beta}=:c_1(n)
 \eeqnn
and
 \beqnn
(x+y^\beta)^{-1} y^{\beta-1+\theta_2}x^{\kappa_2}
 \ar=\ar
(x+y^\beta)^{-(1-\delta_2)}(x+y^\beta)^{-\delta_2} y^{\beta-1+\theta_2}x^{\kappa_2} \cr
 \ar\le\ar
y^{\beta-1+\theta_2-(1-\delta_2)\beta}x^{\kappa_2-\delta_2}
\le n^{\theta_2-1+\delta_2\beta+\kappa_2-\delta_2}
=:c_2(n),
 \eeqnn
which imply
 \beqlb\label{2.20}
(x+y^\beta)^{-1}[\eta_1x^{\theta_1}y^{\kappa_1}+\eta_2\beta y^{\beta-1+\theta_2}x^{\kappa_2}]
\le \eta_1c_1(n)+\eta_2\beta c_2(n),\quad 0<x,y< n.
 \eeqlb
Observe that for all $0<x,y< n$,
 \beqlb\label{2.21}
a_1x^{p_1+1}+b_1\beta y^{q_1+\beta}
\le
a_1 n^{p_1}x+b_1\beta n^{q_1}y^{\beta}
\le
(a_1 n^{p_1}+b_1\beta n^{q_1})(x+y^{\beta})
 \eeqlb
and
 \beqlb\label{2.22}
a_2x^{p_2+2}+b_2\beta^2 y^{q_2+2\beta}
\le
a_2n^{p_2}x^2+b_2\beta^2 n^{q_2}y^{2\beta}
\le
(a_2n^{p_2}+b_2\beta^2 n^{q_2})(x+y^{\beta})^2.
 \eeqlb

Recall \eqref{3.3} and \eqref{3.4}.
By \eqref{3.5}, for $0<x<1$ and $0<y<n$,
 \beqnn
\int_0^1
K_z^1g_n(x,y)\mu_1(\dd z)
 \ar=\ar
\int_0^1
[\ln(1+z/(x+y^\beta))^{-1}+z(x+y^\beta)^{-1}]\mu_1(\dd z) \cr
 \ar\le\ar
\int_0^\infty
[\ln(1+z/(x+y^\beta))^{-1}+z(x+y^\beta)^{-1}]\mu_1(\dd z) \cr
 \ar=\ar
(x+y^\beta)^{-\alpha_1}\int_0^\infty[\ln(1+z)^{-1}+z]\mu_1(\dd z)
=:c_3(x+y^\beta)^{-\alpha_1}.
 \eeqnn
By \eqref{4.1} and \eqref{2.19}, for $0<x<1$ and $0<y<n$
 \beqnn
\int_1^{3n^{\beta}}
K_z^1g_n(x,y)\mu_1(\dd z)
=
\int_1^{3n^{\beta}}z^2\mu_1(\dd z)
\int_0^1(1-u)\partial_{xx} g_n(x+zu,y)\dd u
\le
c_0\int_1^{3n^{\beta}}z^2\mu_1(\dd z).
 \eeqnn
Since $g_n(u,y)=0$ when $u>3n^\beta$,
then $\int_{3n^{\beta}}^\infty
K_z^1g_n(x,y)\mu_1(\dd z)=0$ for all $0<x<1$.
 Thus, it follows that for $0<x<1$ and $0<y<n$,
 \beqnn
x^{p_3+\alpha_1}\int_0^\infty
K_z^1g_n(x,y)\mu_1(\dd z)
\le c_3x^{p_3+\alpha_1}(x+y^\beta)^{-\alpha_1}
+c_0 x^{p_3+\alpha_1}\int_1^{3n^{\beta}}z^2\mu_1(\dd z)
\le c_{31}(n)
 \eeqnn
for some constant $c_{31}(n)>0$.
By \eqref{4.1} and \eqref{2.19} again,
 \beqnn
\int_0^\infty
K_z^1g_n(x,y)\mu_1(\dd z)
 \ar=\ar
\int_0^{3n^{\beta}}z^2\mu_1(\dd z)
\int_0^1(1-u)\partial_{xx} g_n(x+zu,y)\dd u \cr
 \ar\le\ar
c_0\int_0^{3n^{\beta}}z^2\mu_1(\dd z),\qquad 1\le x<n,~0<y<n,
 \eeqnn
and then
 \beqnn
x^{p_3+\alpha_1}\int_0^\infty
K_z^1g_n(x,y)\mu_1(\dd z)\le
c_0n^{p_3+\alpha_1}\int_0^{3n^{\beta}}z^2\mu_1(\dd z)
=:c_{32}(n),\quad 1\le x<n,~0<y<n.
 \eeqnn
Therefore,
 \beqlb\label{2.23}
x^{p_3+\alpha_1}\int_0^\infty
K_z^1g_n(x,y)\mu_1(\dd z)
\le
c_{31}(n)+c_{32}(n)
=:c_3(n),\quad 0<x,y< n.
 \eeqlb
By the same argument, there is a constant $c_4(n)>0$ such that
 \beqlb\label{2.24}
y^{q_3+\alpha_2}\int_0^\infty K_z^2g_n(x,y)\mu_2(\dd z)
\le
c_4(n),\quad 0<x,y<n.
 \eeqlb
Combining \eqref{2.20}--\eqref{2.24} with \eqref{10.1} one has
 \beqlb\label{2.8}
\mathcal{L}g_n(x,y)
 \ar=\ar
(x+y^\beta)^{-1}[\eta_1x^{\theta_1}y^{\kappa_1}+\eta_2\beta y^{\beta-1+\theta_2}x^{\kappa_2}]
+(x+y^\beta)^{-1}[a_1x^{p_1+1}+b_1\beta y^{q_1+\beta}] \cr
 \ar\ar
+(x+y^\beta)^{-2}[a_2x^{p_2+2}+b_2\beta^2 y^{q_2+2\beta}]
-(x+y^\beta)^{-1}\beta(\beta-1) b_2  y^{q_2+\beta} \cr
 \ar\ar
+a_3x^{p_3+\alpha_1}\int_0^\infty
K_z^1g(x,y)\mu_1(\dd z)
+b_3y^{q_3+\alpha_2}\int_0^\infty K_z^2g(x,y)\mu_2(\dd z) \cr
 \ar\le\ar
\eta_1c_1(n)+\eta_2\beta c_2(n)+a_1n^{p_1}+b_1\beta n^{q_1}
+a_2n^{p_2}+b_2\beta^2n^{q_2}\cr
 \ar\ar
+[c_3(n)+c_4(n)]=:C(n),\qquad 0<x,y\le n.
 \eeqlb

Define stopping times  $\tau_{k^{-1}}:=\inf\{t\ge0:X_t+Y_t^\beta\le k^{-1}\}$
and $\tau_0:=\inf\{t\ge0:X_t=Y_t=0\}$.
Let $\gamma_{m,n,k}:=\tau_{m,n}\wedge\tau_{k^{-1}}$ with
$\tau_{m,n}:=\tau_{m^{-1}}^-(X)\wedge\tau_{m^{-1}}^-(Y)\wedge\tau_n^+$
and $\tau_n^+:=\tau_n^+(X)\wedge\tau_n^+(Y)$.
In the following, we estimate $g_n(X_{t\wedge\gamma_{m,n,k}},Y_{t\wedge\gamma_{m,n,k}})$.
By \eqref{a1.1} and It\^o's formula,
 \beqlb\label{2.8b}
g_n(X_{t\wedge\gamma_{m,n,k}},Y_{t\wedge\gamma_{m,n,k}})
=g_n(X_0,Y_0)+\int_0^{{t\wedge\gamma_{m,n,k}}}
\mathcal{L}g_n(X_s,Y_s)\dd s+M(t\wedge\gamma_{m,n,k}),
 \eeqlb
where
 \beqnn
M(t)
 \ar:=\ar
-\int_0^t(X_{s-}+Y_{s-}^\beta)^{-1}\sqrt{2a_2 X_{s-}^{p_2+2}}\dd B_1(s)-\beta\int_0^t
\frac{Y_{s-}^{\beta-1}}{X_{s-}+Y_{s-}^\beta}\sqrt{2b_2 Y_{s-}^{q_2+2}}\dd B_2(s) \cr
 \ar\ar
-\int_0^t\int_0^\infty \int_0^{a_3X_{s-}^{p_3+\alpha_1}}
\ln(1+z/(X_{s-}+Y_{s-}^\beta))\tilde{N}_1(\dd s,\dd z,\dd u) \cr
 \ar\ar
-\int_0^t\int_0^\infty \int_0^{b_3 Y_{s-}^{q_3+\alpha_2}}
\ln\Big[1+\frac{(Y_{s-}+z)^\beta-(Y_{s-})^\beta}{X_{s-}+(Y_{s-})^\beta}\Big]\tilde{N}_2(\dd s,\dd z,\dd u)
=:\sum_{k=1}^4M_{k}(t).
 \eeqnn
In the following, we first show that for each $k=1,2,3,4$, $t\mapsto M_k({t\wedge\gamma_{m,n,k}})$ is a martingale.
For $s<\gamma_n$,
 \beqnn
2a_2 (X_s+Y_s^\beta)^{-2}X_s^{p_2+2}\le
2a_2X_s^{p_2}\le 2a_2n^{p_2},
 \eeqnn
which implies
\beqnn
\int_0^{t\wedge\gamma_{m,n,k}}
\Big[(X_s+Y_s^\beta)^{-1}\sqrt{2a_2 X_s^{p_2+2}}\Big]^2\dd s
\le 2a_2n^{p_2}t.
 \eeqnn
Thus, $t\mapsto M_1({t\wedge\gamma_{m,n,k}})$ is a martingale.
Similarly, $t\mapsto M_2({t\wedge\gamma_{m,n,k}})$ is also a martingale.
Since $\ln(1+z)\le z$ for all $z\ge 0$,
then for $s<\gamma_n$,
 \beqnn
a_3X_s^{p_3+\alpha_1}|\ln(1+z/(X_s+Y_s^\beta))|^2\le a_3z^2X_s^{p_3+\alpha_1}/(X_s+Y_s^\beta)^2
\le z^2a_3n^{p_3+\alpha_1}[n^{-1}+n^{-\beta}]^{-2}
 \eeqnn
and
 \beqnn
a_3X_s^{p_3+\alpha_1}|\ln(1+z/(X_s+Y_s^\beta))|
\le za_3X_s^{p_3+\alpha_1}/(X_s+Y_s^\beta)
\le za_3n^{p_3+\alpha_1}[n^{-1}+n^{-\beta}]^{-1},
 \eeqnn
which imply
 \beqnn
 \ar\ar
\int_0^{t\wedge\gamma_{m,n,k}}\int_0^1 \int_0^{a_3X_s^{p_3+\alpha_1}}
|\ln(1+z/(X_s+Y_s^\beta))|^2\dd s\mu_1(\dd z)\dd u \cr
 \ar\ar\qquad\le
a_3n^{p_3+\alpha_1}[n^{-1}+n^{-\beta}]^{-2}t\int_0^1 z^2\mu_1(\dd z)
 \eeqnn
and
 \beqnn
  \ar\ar
\int_0^{t\wedge\gamma_{m,n,k}}\int_1^\infty \int_0^{a_3X_s^{p_3+\alpha_1}}
|\ln(1+z/(X_s+Y_s^\beta))|\dd s\mu_1(\dd z)\dd u \cr
 \ar\ar\qquad\le
a_3n^{p_3+\alpha_1}[n^{-1}+n^{-\beta}]^{-1}t\int_1^\infty z\mu_1(\dd z).
 \eeqnn
Thus, $t\mapsto M_3({t\wedge\gamma_{m,n,k}})$ is a martingale.
By a similar argument, $t\mapsto M_4({t\wedge\gamma_{m,n,k}})$ is also a martingale.

Combining \eqref{2.8b} with \eqref{2.8} one obtains
 \beqnn
\mbf{E}\big[g_n(X_{t\wedge\gamma_{m,n,k}},Y_{t\wedge\gamma_{m,n,k}})\big]
\le
g_n(X_0,Y_0)+tC(n).
 \eeqnn
It follows from Fatou's lemma that
 \beqnn
 \ar\ar
[\ln (n+n^{\beta})+\ln k]\mbf{P}\{\tau_{k^{-1}}\le t\wedge\tau_0^-\wedge\tau_n^+\}
\le
\mbf{E}\big[g_n(X_{t\wedge\tau_0^-\wedge\tau_n^+\wedge\tau_{k^{-1}}},
Y_{t\wedge\tau_0^-\wedge\tau_n^+\wedge\tau_{k^{-1}}})\big] \cr
 \ar\le\ar
\mbf{E}\big[\lim_{m\to\infty}g_n(X_{t\wedge\gamma_{m,n,k}},
Y_{t\wedge\gamma_{m,n,k}})\big]
\le
\liminf_{m\to\infty}\mbf{E}\big[g_n(X_{t\wedge\gamma_{m,n,k}},
Y_{t\wedge\gamma_{m,n,k}})\big]
\le
g_n(X_0,Y_0)+tC(n).
 \eeqnn
Then
 \beqnn
\lim_{k\to\infty}\mbf{P}\{\tau_{k^{-1}}\le t\wedge\tau_0^-\wedge\tau_n^+\}=0.
 \eeqnn
Taking $n,t\to\infty$ and using Lemma \ref{t10.1}(i), we have $\mbf{P}\{\tau_0\le\tau_0^-, \tau_0<\infty\}=0$ since
	$\tau_0=\lim_{k\to\infty}\tau_{k^{-1}}.$
Notice that $\{\tau_0<\infty\}=\{\tau_0^-(X)=\tau_0^-(Y)<\infty\}$, and
by Definition \ref{def} we have $\tau_0=\tau_0^-$ when $\tau_0<\infty$.
Then
\[\{\tau_0\le\tau_0^-, \tau_0<\infty\}=\{ \tau_0<\infty\}=\{\tau_0^-(X)=\tau_0^-(Y)<\infty\},\]
which implies $ \mbf{P}\{\tau_0^-(X)=\tau_0^-(Y)<\infty\}=0$.
\qed

We first prove the second result of Theorem \ref{t0.01}  by applying Proposition \ref{t3.1b} to the function $g_n(x,y)=g_n(x)=1-\ln (x/n)$ for all $0<x<n$ and $y>0$.

\blemma\label{t3.3}
$\mbf{P}\{\tau_0^-(X)<\infty\}=0$
if $\theta_1\ge1$, and $\mbf{P}\{\tau_0^-(Y)<\infty\}=0$
if $\theta_2\ge1$.
\elemma
\proof
We first assume that $\theta_1\ge1$.
For $n\ge1$ let $g_n\in C^2((0,\infty))$ be a nonnegative function
defined as $g_n(x)=1-\ln (x/n)$ for $0<x<n$
and $g'_n(x)=0$ for all $x\ge n+1$.
Recall \eqref{3.2}.
Then $K_zg_n(x)=0$ for all $x\ge n+1$ and $z>0$.
By \eqref{4.1}, for $0<x\le 1/2$,
 \beqnn
\int_0^{1/2} |K_zg_n(x)|\mu_1(\dd z)
 \ar=\ar
-\int_0^{1/2}[\ln (x+z)-\ln x-zx^{-1}] \mu_1(\dd z) \cr
 \ar=\ar
\int_0^{1/2}\mu_1(\dd z)\int_0^1 z^2 (x+zv)^{-2}(1-v)\dd v \cr
 \ar\le\ar
\int_0^{\infty}\mu_1(\dd z)\int_0^1 z^2 (x+zv)^{-2}(1-v)\dd v.
 \eeqnn
Now applying \eqref{3.5} we get
 \beqnn
\int_0^{1/2} |K_zg_n(x)|\mu_1(\dd z)
\le
x^{-\alpha_1}
\int_0^{\infty}\mu_1(\dd z)\int_0^1 z^2 (1+zv)^{-2}(1-v)\dd v
=:c_{1} x^{-\alpha_1},~~ 0<x\le 1/2.
 \eeqnn
By \eqref{3.2},
 \beqnn
 \ar\ar
\int_{1/2}^\infty K_zg_n(x)\mu_1(\dd z)
 =
\int_{1/2}^{n+1} [g_n(x+z)-g_n(x)-zg_n'(x)]\mu_1(\dd z) \cr
 \ar\le\ar
\sup_{u\ge1/2}|g_n(u)|\int_{1/2}^{n+1} \mu_1(\dd z)
+x^{-1} \int_{1/2}^{n+1} z\mu_1(\dd z)
\le c_{2,n}+c_{2,n}x^{-1},\quad 0<x\le 1/2
 \eeqnn
for some constant $c_{2,n}>0$.
 Thus, it follows that
 \beqlb\label{2.25}
\int_0^\infty |K_zg_n(x)|\mu_1(\dd z)
\le
c_{1} x^{-\alpha_1}+c_{2,n}x^{-1}+c_{2,n},\quad 0<x\le 1/2.
 \eeqlb
Using \eqref{4.1} again,
 \beqnn
\int_0^\infty |K_zg_n(x)|\mu_1(\dd z)
 \ar=\ar
\int_0^{n+1} |K_zg_n(x)|\mu_1(\dd z)
=
\int_0^{n+1}\mu_1(\dd z)\int_0^1 z^2|g''_n(x+zv)|(1-v)\dd v \cr
 \ar\le\ar
\sup_{v\ge1/2}|g''_n(v)|\int_0^{n+1} z^2\mu_1(\dd z)
=:c_{3,n},\qquad x>1/2.
 \eeqnn
Combining this with \eqref{2.25} we have
 \beqnn
\int_0^\infty K_zg_n(x)\mu_1(\dd z)
\le c_{1,n} x^{-\alpha_1}+c_{2,n}x^{-1}+c_{2,n}+c_{3,n},\qquad x>0.
 \eeqnn
Let $g_n(x,y)=g_n(x)$ for $x,y>0$. Since $\theta_1\ge1$, then by \eqref{10.1},
 \beqnn
\mathcal{L}g_n(x,y)
=
\eta_1x^{\theta_1-1}y^{\kappa_1}+G_n(x)
\le
\eta_1n^{\theta_1-1+\kappa_1}+G_n(n)
=:d_n
\le d_n g_n(x,y),\quad 0<x,y\le n,
 \eeqnn
where $G_n(x):=a_1x^{p_1}
+a_2x^{p_2}+a_3[c_{1} x^{-\alpha_1}+ c_{2,n}x^{-1}
+c_{2,n}+c_{3,n}]x^{p_3+\alpha_1}$.
Then by Proposition \ref{t3.1b} and Lemmas \ref{t10.1}(i) and \ref{t3.2}, $\mbf{P}\{\tau_0^-(X)<\infty\}=0$.
By the same argument, we have
$\mbf{P}\{\tau_0^-(Y)<\infty\}=0$ for $\theta_2\ge1$.
This concludes the proof.
\qed

For the proof of the first part of Theorem \ref{t0.01},
we show the following lemma, whose proof is a modification of
that of \cite[Theorem 1.5]{RXYZ19} and we only outline the calculations.

\blemma\label{t3.4}
If $0\le\theta_1\wedge\theta_2<1$,
then
$\mbf{P}\{\tau_0^-<\infty\}>0$.
\elemma
\proof
We only prove the result for $0\le\theta_2<1$.
For simplicity, we assume that $0<X_0,Y_0<1$.
We use Proposition \ref{t0.3} to complete the proof,
which is similar to that of \cite[Theorem 1.5]{RXYZ19}.
Let $0<c_1<c_2<1$ and $\delta>(\theta_1-1)\vee(p+1-\theta_1)\vee1$.
Let $g_0\in C^2((0,1))$
satisfy $g_0(x)=x^{-\delta}$ for $x\in (0,c_1)$
and $g_0(x)=(x-1)^{-2}$ for $x\in(c_2,1)$.
We choose the function $g_0$ so that $g_0,g_0'$ and $g_0''$ are bounded in $[c_1,c_2]$.
For $\lambda_1,\lambda_2>1$ and $0<\rho<1-\theta_2$, define a nonnegative function 
 \beqlb\label{2.0}
g(x,y)
:=
\exp\{-\lambda_1 g_0(x)-\lambda_2 [\tan (y\pi/2) ]^\rho\}1_{\{x,y<1\}},
\qquad
x,y>0,
 \eeqlb
where we only need the properties of the $\tan$ function such
that it is equivalent to $x$ near zero and is infinite at $\pi/2$.
Then by the same arguments as in the proof of \cite[Theorem 1.5]{RXYZ19}, there are constants $\lambda_1,\lambda_2,d>0$ so that
$\mathcal{L}g(x,y)\ge d g(x,y)$ for all $0<x,y<1$,
which can be shown as follows:
(i) For all $\lambda_2>0$ and large enough $\lambda_1>0$, there are constants $d_1=d_1(\lambda_1)>0$ and $d_2:=d_2(\lambda_1)>0$
so that
$\mathcal{L}_1g(x,y)\ge -d_2g(x,y)$
for all $x\ge c_1/2$ and $0<y<1$,
and $\mathcal{L}_1g(x,y)\ge d_1g(x,y)$
for all $0<x<c_1/2$ and $0<y<1$;
(ii) For all $\lambda_1>0$ and large enough $\lambda_2>0$, we have
$\mathcal{L}_2g(x,y)\ge 2d_2g(x,y)$
when $x\ge c_1/2$ and $0<y<1$.
Here
 \beqnn
\mathcal{L}_1g(x,y)
 \ar:=\ar
-\eta_1x^{\theta_1}y^{\kappa_1}g_x'(x,y)-a_1x^{p_1+1}g_x'(x,y) \cr
 \ar\ar
+a_2x^{p_2+2}g_{xx}''(x,y)
+a_3x^{p_3+\alpha_1}\int_0^\infty K_z^1g(x,y)\mu_1(\dd z)
 \eeqnn
and
 \beqnn
\mathcal{L}_2g(x,y)
 \ar:=\ar
-\eta_2y^{\theta_2}x^{\kappa_2}g_y'(x,y)-b_1y^{q_1+1}g_y'(x,y) \cr
 \ar\ar
+b_2y^{q_2+2}g_{yy}''(x,y)
+b_3y^{q_3+\alpha_2}\int_0^\infty K_z^2g(x,y)\mu_2(\dd z).
 \eeqnn
Thus, $\mathcal{L}g(x,y)=\mathcal{L}_1g(x,y)+\mathcal{L}_2g(x,y)
\ge (d_1\wedge d_2)g(x,y)$ for all $x,y>0$.
Applying Proposition \ref{t0.3} one reaches the conclusion.
\qed

Now we are ready to conclude the proof of Theorem \ref{t0.01}.

{\it Proof of Theorem \ref{t0.01}.}
By Lemma \ref{t3.4}, $\mbf{P}\{\tau_0^-<\infty\}>0$ when
$\theta_1\wedge\theta_2<1$.
Then we conclude the conclusion of Lemma \ref{t3.3}.
\qed

\subsection{Proof of Theorem \ref{t0.02}}\label{subsection3.3}

In this subsection, we use Propositions \ref{pro2.2a} and \ref{Pro2.1}
to establish the proof of Theorem \ref{t0.02}.
For $\beta>0$, $0<\delta<\beta^{-1}$ and $0<\rho<1$
we define the key functions $g$ and $\bar{g}$ in Propositions \ref{Pro2.1} by
 \beqlb\label{11.1}
g(x,y):=x^{\beta\delta}y^{-\delta}+y^{\rho},~~\bar{g}(x):=x^{-\delta},\qquad x,y>0.
 \eeqlb
To check condition (ii) in Proposition \ref{Pro2.1}
we first need to compute $g_x',g_y',g_{xx}'',g_{yy}''$
and $\int_0^\infty K_z^1g(x,y)\mu_1(\dd z)$ and $\int_0^\infty K_z^2g(x,y)\mu_2(\dd z)$.

\blemma\label{t4.1}
For $x,y>0$, we have
 \beqlb\label{11.6}
g_x'(x,y)=\beta\delta x^{\beta\delta-1}y^{-\delta},\quad
g_y'(x,y)=-\delta x^{\beta\delta}y^{-\delta-1}+\rho y^{\rho-1}
 \eeqlb
and
 \beqlb\label{11.7}
g_{xx}''(x,y)=-\beta\delta(1-\beta\delta) x^{\beta\delta-2}y^{-\delta},\quad
g_{yy}''(x,y)=\delta(\delta+1) x^{\beta\delta}y^{-\delta-2}-\rho(1-\rho) y^{\rho-2}.
 \eeqlb
Moreover,
 \beqlb\label{11.8a}
\int_0^\infty K_z^1g(x,y)\mu_1(\dd z)
 =
-
\frac{\beta\delta(1-\beta\delta)\Gamma(\alpha_1-\beta\delta)}
{\Gamma(\alpha_1)\Gamma(2-\beta\delta)} x^{\beta\delta-\alpha_1}y^{-\delta}
 \eeqlb
and
 \beqlb\label{11.8b}
\int_0^\infty K_z^2g(x,y)\mu_2(\dd z)
=
\frac{\delta(\delta+1)\Gamma(\alpha_2+\delta)}
{\Gamma(\alpha_2)\Gamma(\delta+2)} x^{\beta\delta}y^{-\delta-\alpha_2}
-\frac{\rho(1-\rho)\Gamma(\alpha_2-\rho)}
{\Gamma(\alpha_2)\Gamma(2-\rho)}y^{\rho-\alpha_2}.
 \eeqlb
\elemma
\proof
The assertions \eqref{11.6} and \eqref{11.7} can be calculated directly
by the definition of $g$ in \eqref{11.1}.
 Using \eqref{3.5} three times and Lemma \ref{t2.0},
 \beqnn
 \ar\ar
\int_0^\infty K_z^1g(x,y)\mu_1(\dd z)
 =
y^{-\delta}\int_0^\infty[(x+z)^{\beta\delta}-x^{\beta\delta}
-\beta\delta zx^{\beta\delta-1}]\mu_1(\dd z) \cr
 \ar=\ar
y^{-\delta}x^{\beta\delta-\alpha_1}\int_0^\infty[(1+z)^{\beta\delta}-1
-\beta\delta z]\mu_1(\dd z)
=
-\frac{\beta\delta(1-\beta\delta)\Gamma(\alpha_1-\beta\delta)}
{\Gamma(\alpha_1)\Gamma(2-\beta\delta)} x^{\beta\delta-\alpha_1}y^{-\delta}
 \eeqnn
and
 \beqnn
 \ar\ar
\int_0^\infty K_z^2g(x,y)\mu_2(\dd z) \cr
 \ar=\ar
x^{\beta\delta}\int_0^\infty[(y+z)^{-\delta}-y^{-\delta}
+\delta zy^{-\delta-1}]\mu_2(\dd z)
+\int_0^\infty[(y+z)^{\rho}-y^{\rho}
-\rho zy^{\rho-1}]\mu_2(\dd z) \cr
 \ar=\ar
x^{\beta\delta}y^{-\delta-\alpha_2}\int_0^\infty[(1+z)^{-\delta}-1
+\delta z]\mu_2(\dd z)
+y^{\rho-\alpha_2}\int_0^\infty[(1+z)^{\rho}-1
-\rho z]\mu_2(\dd z)  \cr
 \ar=\ar
\frac{\delta(\delta+1)\Gamma(\alpha_2+\delta)}
{\Gamma(\alpha_2)\Gamma(\delta+2)} x^{\beta\delta}y^{-\delta-\alpha_2}
-\frac{\rho(1-\rho)\Gamma(\alpha_2-\rho)}
{\Gamma(\alpha_2)\Gamma(2-\rho)}y^{\rho-\alpha_2}.
 \eeqnn
Then \eqref{11.8a} and \eqref{11.8b} follow.
\qed

The key to verifying condition (ii) of Proposition \ref{Pro2.1}
is to verify \eqref{10.9}.
The key to proving \eqref{10.9} is to show that $\tilde{H}_\sigma(x,y)>0$ (to be defined in the following lemma) in a specified region $x,y$
under the assumptions of Theorem \ref{t0.02}.
In the proof of Theorem \ref{t0.02},
\eqref{4.10} and \eqref{4.11} together
imply \eqref{10.9}.

\blemma\label{t2.5}
Suppose that the assumptions in Theorem \ref{t0.02} hold.
Then there are constants $\beta, z_*>0$, $0<\rho,\sigma_0,\varepsilon_0<1$
and $0<\delta<\beta^{-1}$ so that
for all $0<\sigma<\sigma_0$, $0<x,y\le \varepsilon_0$ and $yx^{-\beta}\ge z_*$ we have
 \beqnn
\tilde{H}_\sigma(x,y)
 \ar:=\ar
\delta x^{\beta\delta}y^{-\delta}[\beta a(1-\sigma)x^{p}
-(1+\sigma)by^{q}+\beta\eta_1x^{\theta_1-1}y^{\kappa_1}
-\eta_2y^{\theta_2-1}x^{\kappa_2}] \cr
 \ar\ar
+b\rho(1-\sigma)y^{\rho+q}+\rho\eta_2y^{\rho+\theta_2-1}x^{\kappa_2}
>0.
 \eeqnn
\elemma
\proof
Let $u=x^{-\beta}y $, $0<\delta<\beta^{-1}$ and $0<\sigma<1/2$  in the following.

(i) Suppose that condition (a) in Theorem \ref{t0.02} holds.
 Since $\theta_1\ge1$ and $0\le\theta_2<1$, then $q+1-\theta_2>0$,
$0\le\theta_1-1<\frac{\kappa_2(q-\kappa_1)}{q+1-\theta_2}$ and
$q-\kappa_1>0$. Then $(\theta_1-1)(q+1-\theta_2)-\kappa_2(q-\kappa_1)<0$.
Moreover,
 \beqnn
(\theta_1-1)(\kappa_1+1-\theta_2)-(\kappa_2+1-\theta_1)(q-\kappa_1)=
(\theta_1-1)(q+1-\theta_2)-\kappa_2(q-\kappa_1)<0,
 \eeqnn
which gives $\frac{\theta_1-1}{q-\kappa_1}<\frac{\kappa_2+1-\theta_1}{\kappa_1+1-\theta_2}$.
Thus, there exists a constant $\beta>0$ such that
 \beqnn
\frac{\theta_1-1}{q-\kappa_1}<\beta<\frac{\kappa_2+1-\theta_1}{\kappa_1+1-\theta_2},
 \eeqnn
which implies
 \beqlb\label{11.9a}
\kappa_2+\beta (\theta_2-1)>\theta_1-1+\beta \kappa_1,\qquad
\beta q>\theta_1-1+\beta \kappa_1.
 \eeqlb
Moreover, there are constants $0<\rho,\rho_1<1$ such that
 \beqlb\label{11.9}
\beta q/(1+\rho_1)>\theta_1-1+\beta \kappa_1
 \eeqlb
and
 \beqlb\label{4.2b}
\rho_1\delta>\rho.
 \eeqlb
Observe that there is a constant $z_*>0$
so that for $0<z_*\le u\le y^{-\rho_1}$, we have $y\le x^{\beta/(1+\rho_1)}$
and then for all $0<\sigma<1/2$
 \beqlb\label{4.3}
\tilde{H}_\sigma(x,y)
 \ar\ge\ar
\delta u^{-\delta}[\beta\eta_1x^{\theta_1-1}y^{\kappa_1}
-(1+\sigma)by^{q}
-\eta_2y^{\theta_2-1}x^{\kappa_2}] \cr
 \ar=\ar
\delta u^{-\delta}[\beta\eta_1u^{\kappa_1}x^{\theta_1-1+\beta\kappa_1}
-(1+\sigma)by^{q}
-\eta_2u^{\theta_2-1}x^{\kappa_2+\beta(\theta_2-1)}] \cr
 \ar\ge\ar
\delta u^{-\delta}[\beta\eta_1z_*^{\kappa_1}x^{\theta_1-1+\beta\kappa_1}
-2bx^{\beta q/(1+\rho_1)}
-\eta_2z_*^{\theta_2-1}x^{\kappa_2+\beta(\theta_2-1)}]\cr
 \ar\ge\ar
\delta
u^{-\delta}x^{\theta_1-1+\beta\kappa_1} G_1(x),\qquad 0<z_*\le u\le y^{-\rho_1}.
 \eeqlb
with
 \beqnn
G_1(x):=\beta\eta_1z_*^{\kappa_1}
-2bx^{\frac{\beta q}{(1+\rho_1)}-(\theta_1-1+\beta\kappa_1)}
-\eta_2z_*^{\theta_2-1}x^{\kappa_2+\beta(\theta_2-1)-(\theta_1-1+\beta\kappa_1)}.
 \eeqnn
In view of \eqref{11.9a} and \eqref{11.9},
there is a constant $0<\varepsilon_0<1$ such that $G_1(\varepsilon_0)>0$ and then
 \beqlb\label{4.4}
G_1(x)\ge G_1(\varepsilon_0)>0,\qquad 0<x<\varepsilon_0.
 \eeqlb
Combining \eqref{4.3} and \eqref{4.4} we obtain
 \beqlb\label{2.9b}
\inf_{0<\sigma<1/2}\tilde{H}_\sigma(x,y)>0, \qquad
0<x<\varepsilon_0,~0<z_*\le u\le y^{-\rho_1}.
 \eeqlb
By the definition of $\tilde{H}_\sigma$, for all $0<\sigma<1/2$
and $u> y^{-\rho_1}$,
 \beqlb\label{2.9a}
\tilde{H}_\sigma(x,y)
 \ar\ge\ar
y^{\rho+q}[b\rho(1-\sigma)-\delta (1+\sigma)u^{-\delta}y^{-\rho}]
+\eta_2y^{\rho+\theta_2-1}x^{\kappa_2}[\rho
-\delta u^{-\delta}y^{-\rho}] \cr
 \ar\ge\ar
2^{-1}y^{\rho+q}[b\rho-3\delta y^{\rho_1\delta-\rho}]
+\eta_2y^{\rho+\theta_2-1}x^{\kappa_2}[\rho
-\delta y^{\rho_1\delta-\rho}]=:G_2(y).
 \eeqlb
By \eqref{4.2b}, there is a constant $\varepsilon_1\in(0,\varepsilon_0)$
such that
$G_2(y)>0$ for all $0<y\le \varepsilon_1$.
Combining this with \eqref{2.9a} one obtains
 \beqlb\label{2.9}
\inf_{\sigma\in(0,1/2)}\tilde{H}_\sigma(x,y)>0,\qquad u> y^{-\rho_1},~0<y\le \varepsilon_1.
 \eeqlb
Combining \eqref{2.9} with \eqref{2.9b} one has $\tilde{H}_\sigma(x,y)>0$
for all $0<\sigma<1/2$, $0<x,y\le \varepsilon_1$, and $u\ge z_*$.

(ii) Suppose that condition (b) in Theorem \ref{t0.02} holds.
Then $1-\theta_2>0$, $q+1-\theta_2>0$ and
 \beqnn
p(1-\theta_2)-q(\kappa_2-p)
=
p(q+1-\theta_2)-\kappa_2q<0,
 \eeqnn
which implies
$p/q<(\kappa_2-p)/(1-\theta_2)$.
Thus,
there is a constant $\beta>0$ so that
$p/q<\beta<(\kappa_2-p)/(1-\theta_2)$,
which implies
 \beqlb\label{4.2aa}
\kappa_2+\beta (\theta_2-1)>p,\qquad
\beta q>p.
 \eeqlb
Moreover, there are constants $\rho,\rho_1>0$ so that
 \beqlb\label{4.2a}
\beta q/(1+\rho_1)>p
 \eeqlb
and \eqref{4.2b} holds.
For $u\le y^{-\rho_1}$ we have $y\le x^{\beta/(1+\rho_1)}$ and then
for $0<z_*<1$,
 \beqnn
\tilde{H}_\sigma(x,y)
 \ar\ge\ar
\delta u^{-\delta}[\beta a(1-\sigma)x^{p}-(1+\sigma)by^{q}-\eta_2y^{\theta_2-1}x^{\kappa_2}] \cr
 \ar\ge\ar
2^{-1}\delta u^{-\delta}[\beta ax^{p}-3bx^{q\beta/(1+\rho_1)}
-\eta_2u^{\theta_2-1}x^{\kappa_2+\beta(\theta_2-1)}] \cr
 \ar\ge\ar
2^{-1}\delta u^{-\delta}x^p\big[\beta a-3bx^{q\beta/(1+\rho_1)-p}
-\eta_2u_*^{\theta_2-1}x^{\kappa_2+\beta(\theta_2-1)-p}\big]
=:2^{-1}\delta u^{-\delta}x^p G_3(x).
\eeqnn
By \eqref{4.2aa} and \eqref{4.2a}, there is a constant $0<\varepsilon_0<1$
such that
 \beqnn
G_3(x)\ge G_3(\varepsilon_0)>0,\qquad 0<x<\varepsilon_0.
 \eeqnn
Then combine this with \eqref{2.9}, $\tilde{H}_\sigma(x,y)>0$ for some constant $\varepsilon_1\in(0,\varepsilon_0)$
when $0<x,y<\varepsilon_1$, $0<\sigma\le1/2$ and $u> y^{-\rho_1}$.

(iii) Suppose that condition (c) in Theorem \ref{t0.02} holds.
Taking $b/a<\beta<\kappa_2/(1-\theta_2)$,
we have $\kappa_2+\beta(\theta_2-1)>0$ and
$\beta a(1-\sigma)>(1+\sigma)b$  for all $0<\sigma\le \sigma_0$ and some $0<\sigma_0<1/2$.
Moreover,
there are constants $z_*,\varepsilon_0\in(0,1)$
so that for all $0<x<\varepsilon_0$ and $u\ge z_*$,
\beqnn
\tilde{H}_\sigma(x,y)
 \ar\ge\ar
\delta u^{-\delta}\big[\beta a(1-\sigma)
-(1+\sigma)b
-\eta_2u^{\theta_2-1}x^{\kappa_2+\beta(\theta_2-1)}\big] \cr
 \ar\ge\ar
\delta u^{-\delta}\big[\beta a(1-\sigma)
-(1+\sigma)b
-\eta_2z_*^{\theta_2-1}\varepsilon_0^{\kappa_2+\beta(\theta_2-1)}\big]>0.
 \eeqnn

(iv) Suppose that condition (d) in Theorem \ref{t0.02} holds. Taking  $b/\eta_1<\beta<\kappa_2/(\kappa_1+1-\theta_2)$,
we have $\beta\eta_1>b$ and $\beta\kappa_1<\kappa_2+\beta(\theta_2-1)$.
It follows that there is a constant $\sigma_0\in(0,1/2)$ such that $\beta\eta_1-(1+\sigma)b>0$
for all $0<\sigma<\sigma_0$.
Moreover, there are  constants $\varepsilon_0\in(0,1/2)$ and $z_*>0$
so that for all $0<\sigma<\sigma_0$, $0<x<\varepsilon_0$ and $u\ge  z_*$, we have
 \beqnn
\tilde{H}_\sigma(x,y)
 \ar\ge\ar
\delta  u^{-\delta}\big[u^{\kappa_1}[\beta\eta_1-(1+\sigma)b]x^{\beta\kappa_1}
-\eta_2u^{\theta_2-1}x^{\kappa_2+\beta(\theta_2-1)}\ \big]\cr
 \ar\ge\ar
\delta  u^{-\delta} x^{\beta\kappa_1}\big[u^{\kappa_1}[\beta\eta_1-(1+\sigma)b]
-\eta_2u^{\theta_2-1}x^{\kappa_2+\beta(\theta_2-1)-\beta\kappa_1} \big] \cr
 \ar\ge\ar
\delta u^{-\delta} x^{\beta\kappa_1}\big[z_*^{\kappa_1}[\beta\eta_1-(1+\sigma)b]
-\eta_2z_*^{\theta_2-1}\varepsilon_0^{\kappa_2+\beta(\theta_2-1)-\beta\kappa_1} \big]>0.
 \eeqnn
This completes the proof.
\qed

Now we are ready to conclude the proof of Theorem \ref{t0.02}.

{\it Proof of Theorem \ref{t0.02}.}
$\mbf{P}\{\tau_0^-(Y)<\infty\}>0$ follows from Lemmas \ref{t3.3} and \ref{t3.4}.
We use Proposition \ref{Pro2.1} to prove the rest of the assertion
with $g,\bar{g}$ defined by \eqref{11.1}.
In the following, we verify the condition \eqref{10.9}.
Let
 \beqnn
H(x,y):=\delta x^{\beta\delta}y^{-\delta}H_1(x,y)+\rho y^{\rho}H_2(x,y)
 \eeqnn
with
 \beqnn
H_1(x,y)
 \ar:=\ar
a_1\beta  x^{p_1}
+\beta (1-\beta\delta)a_2x^{p_2}
+
\frac{\beta (1-\beta\delta)\Gamma(\alpha_1-\beta\delta)}
{\Gamma(\alpha_1)\Gamma(2-\beta\delta)}a_3x^{p_3} \cr
 \ar\ar
-  b_1 y^{q_1}
- (\delta+1)b_2y^{q_2}
-\frac{\delta(\delta+1)\Gamma(\alpha_2+\delta)}
{\Gamma(\alpha_2)\Gamma(2+\delta)}b_3y^{q_3}
+\beta\delta\eta_1x^{\theta_1-1}y^{\kappa_1}
-\delta\eta_2y^{\theta_2-1}x^{\kappa_2}
 \eeqnn
and
 \beqnn
H_2(x,y):= b_1 y^{q_1}+(1-\rho)y^{q_2}
+\frac{(1-\rho)\Gamma(\alpha_2-\rho)}
{\Gamma(\alpha_2)\Gamma(2-\rho)}b_3y^{q_3}
+\eta_2y^{\theta_2-1}x^{\kappa_2}.
 \eeqnn
Let $\tilde{H}_\sigma(x,y)$ be determined in Lemma \ref{t2.5}.
By Lemma \ref{t2.5},
there are constants $\beta,\rho,\delta,\sigma,\varepsilon_0,u_*>0$ such that
 \beqlb\label{4.10}
H(x,y)\ge \tilde{H}_\sigma(x,y)\ge0,\qquad  0<x,y\le \varepsilon_0,~~yx^{-\beta}\ge u_*.
 \eeqlb
By \eqref{10.1} and Lemma \ref{t4.1},
 \beqlb\label{4.11}
\mathcal{L}g(x,y)
\le-H(x,y),\qquad 0<x,y\le \varepsilon_0,~x^{-\beta}y\ge u_*,
 \eeqlb
which gives condition \eqref{10.9}.
 Condition \eqref{10.9b} is
obtained by the definition of function
$g$ in \eqref{11.1} and \eqref{10.1} again.
 Condition (iii) of Proposition \ref{Pro2.1}
is met by Lemmas \ref{t10.1} and \ref{t3.3}.
Let $u_0>u_*$ be fixed and let $\tilde{\varepsilon},\tilde{\varepsilon}_0\in(0,\varepsilon_0)$
be small enough so that
$g(X_0,Y_0)<(1-\tilde{\varepsilon})u_*^{-\delta}$ for $X_0= \varepsilon^{1+\beta^{-1}}$ and $Y_0=u_0\varepsilon^{\beta+1}$ with $0<\varepsilon<\tilde{\varepsilon}_0$.
Now, using Proposition \ref{Pro2.1},
$\mbf{P}\{\tau_0^-(Y)<\infty\}<1$ for $X_0= \varepsilon^{1+\beta^{-1}},Y_0=u_0\varepsilon^{\beta+1}$, $0<\varepsilon<\varepsilon_1$ and some $0<\varepsilon_1<\tilde{\varepsilon}_0$.

By the comparison theorem (Proposition \ref{t2.1}) and Proposition \ref{Pro2.1},
$\mbf{P}\{\tau_0^-(Y)=\infty\}>0$
for $X_0\le \varepsilon^{1+\beta^{-1}}$ and
$Y_0\ge u_0\varepsilon^{\beta+1}$.
For general $X_0,Y_0>0$ let $D:=\{(x,y):x\le\varepsilon^{1+\beta^{-1}},
y\ge u_0\varepsilon^{\beta+1}\}$ and $\tau_D:=\inf\{t\ge0:(X_t,Y_t)\in D\}$.
It follows from
Proposition \ref{t2.10} that $\mbf{P}\{\tau_D<\infty\}>0$.
Applying the strong Markov property,
 \beqnn
\mbf{P}\{\tau_0^-(Y)=\infty\}
 \ar\ge\ar
\mbf{P}\{\tau_0^-(Y)\circ\hat{\theta}(\tau_D)=\infty,\tau_D<\infty\} \cr
 \ar=\ar
\mbf{E}\Big[\mbf{E}\big[1_{\{\tau_0^-(Y)\circ\hat{\theta}(\tau_D)=\infty\}}
1_{\{\tau_D<\infty\}}|\mathcal{F}_{\tau_D}\big]\Big] \cr
 \ar=\ar
\mbf{E}\Big[\mbf{P}\big\{\tau_0^-(Y)=\infty
|X_0\le \varepsilon^{1+\beta^{-1}},Y_0\ge u_0\varepsilon^{\beta+1}\big\}
1_{\{\tau_D<\infty\}}\Big] \cr
 \ar=\ar
\mbf{P}\big\{\tau_0^-(Y)=\infty
|X_0\le \varepsilon^{1+\beta^{-1}},Y_0\ge u_0\varepsilon^{\beta+1}\big\}
\cdot \mbf{P}\{\tau_D<\infty\}>0,
 \eeqnn
where $\hat{\theta}(t)$ denotes the usual shift operator.
This ends the proof.
\qed

\subsection{Proof of Theorem \ref{t0.03}}

In this section, we prove Theorem \ref{t0.03} by using Corollary \ref{cor2.1}
with the function $g$ given by
 \beqnn
g(x,y):=g(yx^{-\beta}),~~g(x):=\e^{-\lambda x^r},\qquad x,y>0,
 \eeqnn
for $0<r<1$ and $\beta,\lambda>0$.
To verify \eqref{2.4a} of Corollary \ref{cor2.1}, we first
 estimate $g_x',g_{xx}'',g_y',g_{y}''$
and $\int_0^\infty K_z^1g(x,y)\mu_1(\dd z)$ and $\int_0^\infty K_z^2g(x,y)\mu_2(\dd z)$.

\blemma\label{t2.6}
For $x,y>0$ let $u:=yx^{-\beta}$. Then
 \beqlb\label{11.20}
g'_x(x,y)=\lambda r\beta u^rg(u)x^{-1},\quad g'_y(x,y)=-\lambda ru^rg(u)y^{-1}
 \eeqlb
and
 \beqlb\label{11.20a}
g''_{xx}(x,y)\ge
-\lambda r\beta(r\beta+1) u^rg(u)x^{-2},\quad
g''_{yy}(x,y)\ge\lambda r(1-r) u^rg(u)y^{-2}.
 \eeqlb
Moreover,
 \beqnn
\int_0^\infty K_z^1g(x,y)\mu_1(\dd z)
\ge
-\lambda r u^r x^{-\alpha_1}\frac{\beta(r\beta+1)\Gamma(\alpha_1+r\beta)}
{\Gamma(\alpha_1)\Gamma(2+r\beta)}
 \eeqnn
and
 \beqnn
\int_0^\infty K_z^2g(x,y)\mu_2(\dd z)
\ge
\lambda r u^r
y^{-\alpha_2}\frac{r(1-r)\Gamma(\alpha_2-r)}
{\Gamma(\alpha_2)\Gamma(2-r)}.
 \eeqnn
\elemma
\proof
It is elementary to see that
 \beqnn
g'_x(x,y)=\lambda r\beta u^rg(u)x^{-1},\quad
g''_{xx}(x,y)=\lambda^2 r^2\beta^2 u^{2r}g(u)x^{-2}
-\lambda r\beta(r\beta+1) u^rg(u)x^{-2}
 \eeqnn
and
 \beqnn
g'_y(x,y)=-\lambda ru^rg(u)y^{-1},\quad
g''_{yy}(x,y)=\lambda^2 r^2u^{2r}g(u)y^{-2}
+\lambda r(1-r) u^rg(u)y^{-2},
 \eeqnn
which give \eqref{11.20} and \eqref{11.20a}.
Since $\e^{x}-1\ge x$ for all $x\in\mbb{R}$, then
 \beqnn
K_z^1g(x,y)
 \ar=\ar
g(u)\big[\e^{-\lambda y^r(x+z)^{-r\beta}+\lambda y^rx^{-r\beta}}-1
-\lambda r\beta u^r x^{-1}z\big] \cr
 \ar\ge\ar
-g(u)\big[\lambda y^r(x+z)^{-r\beta}-\lambda y^rx^{-r\beta}
+\lambda r\beta u^r x^{-1}z\big] \cr
 \ar=\ar
-\lambda u^rg(u)\big[(1+x^{-1}z)^{-r\beta}-1
+r\beta x^{-1}z\big].
 \eeqnn
Since $\mu_1(\dd z)=\frac{\alpha_1(\alpha_1-1)}{\Gamma(\alpha_1)\Gamma(2-\alpha_1)}z^{-1-\alpha_i}\dd z$, replacing the variable $z$ with $zx$, we have
 \beqnn
\int_0^\infty\big[(1+x^{-1}z)^{-r\beta}-1
+r\beta x^{-1}z\big]\mu_1(\dd z)
=x^{-\alpha_1}
\int_0^\infty\big[(1+z)^{-r\beta}-1
+r\beta z\big]\mu_1(\dd z).
 \eeqnn
 Thus, from Lemma \ref{t2.0} it follows that
 \beqnn
[\lambda u^rg(u)] ^{-1} \int_0 ^ \infty K_z^1g(x,y) \mu_1 (\dd z)
 \ge
-x^{-\alpha_1}\frac{r\beta(r\beta+1)\Gamma(\alpha_1+r\beta)}
{\Gamma(\alpha_1)\Gamma(2+r\beta)}.
 \eeqnn
Similarly,
 \beqnn
K_z^2g(x,y)
 \ar=\ar
g(u)\big[\e^{-\lambda (y+z)^rx^{-r\beta}+\lambda y^rx^{-r\beta}}-1
+\lambda r u^r y^{-1}z\big] \cr
 \ar\ge\ar
-g(u)\big[\lambda (y+z)^rx^{-r\beta}-\lambda y^rx^{-r\beta}
-\lambda r u^r y^{-1}z\big] \cr
 \ar=\ar
-\lambda u^rg(u)\big[(1+y^{-1}z)^{r}-1
-r y^{-1}z\big],
 \eeqnn
which implies
\beqnn
 \ar\ar
[\lambda u^rg(u)]^{-1}\int_0^\infty K_z^2g(x,y)\mu_2(\dd z) \cr
 \ar\ge\ar
- y^{-\alpha_2}
\int_0^\infty\big[(1+z)^{r}-1
-r z\big]\mu_2(\dd z)
=y^{-\alpha_2}\frac{r(1-r)\Gamma(\alpha_2-r)}
{\Gamma(\alpha_2)\Gamma(2-r)}.
 \eeqnn
This ends the proof.
\qed

To verify \eqref{2.4a}, the key is to show \eqref{10.11} in the following Lemma \ref{t2.7} under the assumptions of Theorem \ref{t0.03}
(see \eqref{3.34} of the proof of Theorem \ref{t0.03}
in the following).

\blemma\label{t2.7}
Suppose that the assumptions in Theorem \ref{t0.03} hold.
Then there are constants $r,\varepsilon\in(0,1)$, $\beta>0$, and $c_0>0$
such that
 \beqlb\label{10.11}
H(x,y)\ge c_0,\qquad 0<x,y<\varepsilon,
 \eeqlb
where $u:=x^{-\beta}y$ and
 \beqnn
H(x,y):=b(r) u^r y^q+\eta_2u^ry^{\theta_2-1}x^{\kappa_2}
-a(r)\beta u^r  x^p
-\eta_1\beta u^r x^{\theta_1-1}y^{\kappa_1}
 \eeqnn
with
 \beqnn
a(r):=
a_11_{\{p_1=p\}}+a_2(r\beta+1)1_{\{p_2=p\}}
+a_3 \frac{(r\beta+1)\Gamma(\alpha_1+r\beta)}
{\Gamma(\alpha_1)\Gamma(2+r\beta)} 1_{\{p_3=p\}}
 \eeqnn
and
 \beqnn
b(r):=
b_11_{\{q_1=q\}}+b_2(1-r)1_{\{q_2=q\}}
+b_3\frac{(1-r)\Gamma(\alpha_2-r)}
{\Gamma(\alpha_2)\Gamma(2-r)} 1_{\{q_3=q\}}.
 \eeqnn
\elemma
\proof
Let $0<r<1-\theta_2$ and
 \beqnn
 \delta_1:=\frac{1-\theta_2-r}{q+1-\theta_2},~~
\delta_2:=\frac{1-\theta_2}{q+1-\theta_2},~~
\delta_3:=\frac{\kappa_1+1-\theta_2}{q+1-\theta_2},~~
c_i(r):=\delta_i^{-\delta_i}(1-\delta_i)^{\delta_i-1}b(r)^{\delta_i} \eta_2^{1-\delta_i}
 \eeqnn
for $i=1,2,3$ and $q>\kappa_1$.
Then $0<\delta_1,\delta_2,\delta_3<1$ and $q\delta_1+(\theta_2-1)(1-\delta_1)+r=
q\delta_2+(\theta_2-1)(1-\delta_2)
=0$.
By Lemma \ref{t2.3},
 \beqlb\label{11.2}
b(r) u^r y^q+\eta_2u^ry^{\theta_2-1}x^{\kappa_2}
 \ar\ge\ar
\delta_1^{-\delta_1}(1-\delta_1)^{\delta_1-1}
[b(r) u^r y^q]^{\delta_1}\cdot [\eta_2u^ry^{\theta_2-1} x^{\kappa_2}]^{1-\delta_1} \cr
 \ar=\ar
c_1(r)u^ry^{q\delta_1+(\theta_2-1)(1-\delta_1)} x^{\kappa_2(1-\delta_1)}
=
c_1(r)x^{\frac{\kappa_2(r+q)}{q+1-\theta_2}-r\beta}
 \eeqlb
and
 \beqlb\label{11.3}
b(r) y^q+\eta_2y^{\theta_2-1}x^{\kappa_2}
 \ar\ge\ar
\delta_2^{-\delta_2}(1-\delta_2)^{\delta_2-1}
[b(r) y^q]^{\delta_2}\cdot [\eta_2y^{\theta_2-1}x^{\kappa_2}]^{1-\delta_2} \cr
 \ar=\ar
c_2(r)y^{q\delta_2+(\theta_2-1)(1-\delta_2)}x^{\kappa_2(1-\delta_2)}
=
c_2(r)x^{\frac{\kappa_2q}{q+1-\theta_2}}.
 \eeqlb
Similarly, for $q>\kappa_1$,
we have $q\delta_3+(\theta_2-1)(1-\delta_3)=\kappa_1$,
$1-\delta_3=\frac{(q-\kappa_1)}{q+1-\theta_2}$ and then
 \beqlb\label{11.3b}
 \ar\ar
b(r) y^q+\eta_2y^{\theta_2-1}x^{\kappa_2}
\ge
\delta_3^{-\delta_3}(1-\delta_3)^{\delta_3-1}
[b(r) y^q]^{\delta_3}\cdot [\eta_2y^{\theta_2-1}x^{\kappa_2}]^{1-\delta_3} \cr
 \ar=\ar
c_3(r)y^{\kappa_1}x^{\kappa_2(1-\delta_3)}
=
c_3(r)u^{\kappa_1}x^{\kappa_2(1-\delta_3)
+\beta\kappa_1}
c_3(r)u^{\kappa_1} x^{\frac{\kappa_2(q-\kappa_1)}{q+1-\theta_2}+\beta\kappa_1}.
 \eeqlb

In the following, we assume that condition (b) of Theorem \ref{t0.03} is satisfied.
Since $a,b>0$, there is a constant $r_1\in(0,1-\theta_2)$
such that $\inf_{r\in(0,r_1]}[a(r)\wedge b(r)] >0$.
Let $\beta> \kappa_2/(q+1-\theta_2)$
satisfy
 \beqlb\label{11.18}
\frac{\kappa_2(r_1+q)}{q+1-\theta_2}-r_1\beta<0.
 \eeqlb
Under conditions (a) and (b) of Theorem \ref{t0.03},
there is a constant $\varepsilon_1\in(0,1)$ so that
 \beqnn
c_2(r_1)[\beta a(r_1)]^{-1}\varepsilon_1^{\frac{\kappa_2q}{q+1-\theta_2}-p}\ge3,\quad
c_3(r_1)(\eta_1\beta)^{-1}
\varepsilon_1^{\frac{\kappa_2(q-\kappa_1)}{q+1-\theta_2}-(\theta_1-1)}
\ge3.
 \eeqnn
 Thus, it follows from \eqref{11.3} and \eqref{11.3b} that
 \beqlb\label{11.4}
b(r_1) y^q+\eta_2y^{\theta_2-1}x^{\kappa_2}
 \ar\ge\ar
\beta a(r_1)x^p\cdot
c_2(r_1)[\beta a(r_1)]^{-1}x^{\frac{\kappa_2q}{q+1-\theta_2}-p} \cr
 \ar\ge\ar
\beta a(r_1)x^p\cdot
c_2(r_1)[\beta a(r_1)]^{-1}\varepsilon_1^{\frac{\kappa_2q}{q+1-\theta_2}-p}
\ge
3\beta a(r_1)x^p
 \eeqlb
and
 \beqlb\label{11.5}
b(r_1) y^q+\eta_2y^{\theta_2-1}x^{\kappa_2}
 \ar\ge\ar
c_3(r_1)u^{\kappa_1} x^{\frac{\kappa_2(q-\kappa_1)}{q+1-\theta_2}+\beta\kappa_1}
\ge
\eta_1\beta x^{\theta_1-1}y^{\kappa_1}\cdot
c_3(r_1)[\eta_1\beta]^{-1} x^{\frac{\kappa_2(q-\kappa_1)}{q+1-\theta_2}-(\theta_1-1)} \cr
 \ar\ge\ar
\eta_1\beta x^{\theta_1-1}y^{\kappa_1}\cdot
c_3(r_1)[\eta_1\beta]^{-1} \varepsilon_1^{\frac{\kappa_2(q-\kappa_1)}{q+1-\theta_2}-(\theta_1-1)}
\ge
3\eta_1\beta x^{\theta_1-1}y^{\kappa_1}
 \eeqlb
for all $0<x<\varepsilon_1$ when $q>\kappa_1$.
When $q<\kappa_1$,
 \beqnn
b(r_1) y^q+\eta_2y^{\theta_2-1}x^{\kappa_2}
\ge b(r_1) y^q\ge
3\eta_1\beta x^{\theta_1-1}y^{\kappa_1},\quad 0<y<\varepsilon_2
 \eeqnn
for some $0<\varepsilon_2<\varepsilon_1$.
When $q=\kappa_1$ and $\theta_1>1$,
 \beqnn
b(r_1) y^q+\eta_2y^{\theta_2-1}x^{\kappa_2}
\ge b(r_1) y^q\ge
3\eta_1\beta x^{\theta_1-1}y^{\kappa_1},\quad 0<x<\varepsilon_3
 \eeqnn
for some $0<\varepsilon_3<\varepsilon_2$.
Thus, \eqref{11.5} holds for all $0<x,y<\varepsilon_3$.
Now by \eqref{11.2} and \eqref{11.18}--\eqref{11.4},
for $r=r_1$, we obtain
 \beqnn
H(x,y)
 \ar=\ar
3^{-1}[b(r) u^r y^q+\eta_2u^ry^{\theta_2-1}x^{\kappa_2}]
-a(r)\beta u^r  x^p \cr
 \ar\ar
+3^{-1}[b(r) u^r y^q+\eta_2u^ry^{\theta_2-1}x^{\kappa_2}]
-\eta_1\beta u^r x^{\theta_1-1}y^{\kappa_1}\cr
 \ar\ar
+3^{-1}[b(r) u^r y^q+\eta_2u^ry^{\theta_2-1}x^{\kappa_2}]
\ge 3^{-1}c_1(r),\quad 0<x,y<\varepsilon_3.
 \eeqnn
This gives \eqref{10.11}.

Suppose that condition (c) of Theorem \ref{t0.03} holds.
Let $\kappa_2/(1-\theta_2)<\beta<b/a$.
Then $\kappa_2+\beta(\theta_2-1)<0$
and $b>a\beta$.
Thus
$b+\eta_2y^{\theta_2-1}x^{\kappa_2}>
a\beta$.
Moreover, for $0<r<1-\theta_2$, $u^r[b+\eta_2y^{\theta_2-1}x^{\kappa_2}]>a\beta$ when $u\ge1$, and
 \beqnn
u^r[b+\eta_2y^{\theta_2-1}x^{\kappa_2}]
\ge
\eta_2u^{r+\theta_2-1} x^{\kappa_2+\beta(\theta_2-1)}
\ge \eta_2 ,\quad 0<x<1
 \eeqnn
when $0<u<1$.
It follows that
 \beqnn
u^r[b+\eta_2y^{\theta_2-1}x^{\kappa_2}]
\ge (a\beta)\wedge\eta_2 ,\qquad u>0,~0<x<1.
 \eeqnn
By the continuities of $r\mapsto a(r)$ and $r\mapsto b(r)$, there are
constants $r_2\in(0,1-\theta_2)$, $\varepsilon,\varepsilon_4\in(0,1/2)$, and $d_1>0$ so that for all $0<r<r_2$ and $0<x,y<\varepsilon$, we have
 \beqnn
(1-2\varepsilon_4)[b(r) +\eta_2y^{\theta_2-1}x^{\kappa_2}]
>a(r)\beta,\quad \eta_1\beta x^{\theta_1-1}y^{\kappa_1}
\le
\eta_1\beta \varepsilon^{\kappa_1}
\le b(r)\varepsilon_4
 \eeqnn
and
 \beqnn
b(r)u^r+\eta_2u^ry^{\theta_2-1}x^{\kappa_2}
\ge
d_1.
 \eeqnn
Thus
 \beqnn
H(x,y)
 \ar=\ar
(1-2\varepsilon_4)[b(r) u^r +\eta_2u^ry^{\theta_2-1}x^{\kappa_2}]
-a(r)\beta u^r  \cr
 \ar\ar
+\varepsilon_4[b(r) u^r +\eta_2u^ry^{\theta_2-1}x^{\kappa_2}]
-\eta_1\beta u^r x^{\theta_1-1}y^{\kappa_1}\cr
 \ar\ar
+\varepsilon_4[b(r) u^r +\eta_2u^ry^{\theta_2-1}x^{\kappa_2}]
\ge \varepsilon_4 d_1,\quad 0<x,y<\varepsilon,
 \eeqnn
which gives \eqref{10.11}.
\qed

Now we are ready to complete the proof of Theorem \ref{t0.03}.

{\it Proof of Theorem \ref{t0.03}.}
Let $g$ be the function defined at the beginning of this subsection.
It follows from \eqref{10.1} and Lemmas \ref{t2.6}--\ref{t2.7} that
 \beqlb\label{3.34}
\mathcal{L}g(x,y) \ge H(x,y) g(x,y)
\ge c_0g(x,y),\qquad 0<x,y\le \varepsilon
 \eeqlb
for some constants $c_0>0$ and $\varepsilon\in(0,1)$.
We first assume that $0<X_0,Y_0\le\varepsilon_1<\varepsilon$.
Combining this with Lemma \ref{t10.1}(i) and Corollary \ref{cor2.1},
we obtain
 \beqnn
\mbf{P}\{\tau_0^-(X)\wedge\tau_0^-(Y)
\wedge\tau_\varepsilon^+(X)\wedge\tau_\varepsilon^+(Y)<\infty\}
\ge\e^{-\lambda(Y_0X_0^{-\beta})^r}.
 \eeqnn
By Lemma \ref{t3.3},
$\tau_0^-(X)=\infty$ almost surely.
Then letting $\lambda\to0$,
 \beqlb\label{3.21}
\mbf{P}\{\tau_0^-(Y)
\wedge\tau_\varepsilon^+(X)\wedge\tau_\varepsilon^+(Y)<\infty\}=1.
 \eeqlb
By Lemma \ref{t10.1}(ii), we have
\beqnn
\mbf{P}\{\tau_{\varepsilon}^+(X)<\infty\}
=\mbf{P}\big\{\sup_{t\ge0}X_t\ge\varepsilon\big\}
\le C(\varepsilon_1/\varepsilon)^{1/4}\quad\text{and}\quad
\mbf{P}\{\tau_{\varepsilon}^+(Y)<\infty\}
\le C(\varepsilon_1/\varepsilon)^{1/4}
\eeqnn
for some constant $C>0$,
which gives
\beqnn
\mbf{P}\{\tau_{\varepsilon}^+(X)\wedge \tau_{\varepsilon}^+(Y)<\infty\}
\le
\mbf{P}\{\tau_{\varepsilon}^+(X)<\infty\}
+\mbf{P}\{\tau_{\varepsilon}^+(Y)<\infty\}
\le 2C(\varepsilon_1/\varepsilon)^{1/4}.
\eeqnn
Since
 \beqnn
\mbf{P}\{\tau_0^-(Y)<\infty\}
\ge
\mbf{P}\{\tau_0^-(Y)\wedge\tau_{\varepsilon}^+(X)\wedge\tau_{\varepsilon}^+(Y)<\infty\}
-\mbf{P}\{\tau_{\varepsilon}^+(X)\wedge\tau_{\varepsilon}^+(Y)<\infty\},
 \eeqnn
then by \eqref{3.21},
\beqlb\label{proof1.6}
\mbf{P}\{\tau_0^-(Y)=\infty\}
\le
2C(\varepsilon_1/\varepsilon)^{1/4}.
\eeqlb

In the following, we assume that $X_0>\varepsilon_1$ or $Y_0>\varepsilon_1$,
and $\sigma:=\inf\{t\ge0:X_t+Y_t\le \varepsilon_1\}$.
By Lemma \ref{t10.1}(i) and the strong Markov property,
 \beqnn
\mbf{P}\{\tau_0^-(Y)=\infty\}
 \ar=\ar
\mbf{P}\{\tau_0^-(Y)=\infty,\sigma<\infty\}
=\mbf{P}\{\tau_0^-(Y)\circ\hat{\theta}(\sigma)=\infty,\sigma<\infty\} \cr
 \ar=\ar
\mbf{E}\Big[\mbf{E}\big[1_{\{\tau_0^-(Y)\circ\hat{\theta}(\sigma)=\infty\}}
1_{\{\sigma<\infty\}}
|\mathcal{F}_\sigma\big]\Big] \cr
 \ar=\ar
\mbf{E}\Big[\mbf{P}\{\tau_0^-(Y)=\infty|X_0+Y_0\le \varepsilon_1\}
1_{\{\sigma<\infty\}}\Big],
 \eeqnn
where we recall that $\hat{\theta}(t)$ is the shift operator.
Applying \eqref{proof1.6}, we have
 \beqnn
\mbf{P}\{\tau_0^-(Y)=\infty\}
\le
2C(\varepsilon_1/\varepsilon)^{1/4}
\mbf{P}\{\sigma<\infty\}\le
2C(\varepsilon_1/\varepsilon)^{1/4},
 \eeqnn
which implies that \eqref{proof1.6} holds for all $X_0,Y_0>0$. Taking
$\varepsilon_1\to0$ in \eqref{proof1.6} one concludes the proof.
\qed

\section{Appendix}\label{Appendix}
\setcounter{equation}{0}

In this section, we present two lemmas which are used in Section
\ref{proof of Theorems}. We also present the proof for Proposition \ref{t2.1}.
 \blemma\label{t2.3}
For any $u,v\ge0$ and $p,q>1$ with $1/p+1/q=1$, we have
 \beqnn
u+v\ge p^{1/p}q^{1/q} u^{1/p}v^{1/q}.
 \eeqnn
\elemma
\proof
It follows from the Young inequality immediately.
\qed

\blemma\label{t2.0}
Given $\mu(\dd z):=\frac{\alpha(\alpha-1)}{\Gamma(\alpha)\Gamma(2-\alpha)}z^{-1-\alpha}\dd z$ for $1<\alpha<2$,
 for $\beta>0$ we have
 \beqlb\label{2.6}
\int_0^\infty \big[1-(1+z)^{-\beta}\big]z \mu(\dd z)
=\frac{\alpha\beta\Gamma(\alpha+\beta-1)}{\Gamma(\alpha)\Gamma(\beta+1)}
 \eeqlb
and
 \beqnn
\int_0^\infty \big[(1+z)^{-\beta}-1+\beta z\big]\mu(\dd z)
=\frac{\beta(\beta+1)\Gamma(\alpha+\beta)}
{\Gamma(\alpha)\Gamma(\beta+2)}.
 \eeqnn
Moreover, for $0<\beta<\alpha-1$ we have
 \beqnn
\int_0^\infty \big[(1+z)^{\beta}-1\big]z\mu(\dd z)
=\frac{\alpha\beta\Gamma(\alpha-\beta-1)}{\Gamma(\alpha)\Gamma(1-\beta)}
 \eeqnn
and for $0<\beta<1$ we have
 \beqnn
\int_0^\infty \big[(1+z)^{\beta}-1-\beta z\big]\mu(\dd z)
=-\frac{\beta(1-\beta)\Gamma(\alpha-\beta)}
{\Gamma(\alpha)\Gamma(2-\beta)}.
 \eeqnn
\elemma
\proof
Since the proofs are similar, we only state that of \eqref{2.6}.
By integration by parts, we have
 \beqnn
 \ar\ar
\int^\infty_0 [1-(1+z)^{-\beta}]z^{-\alpha} \dd z
=(\alpha-1)^{-1}\beta
\int^\infty_0 (1+z)^{-\beta-1}z^{1-\alpha} \dd z \cr
 \ar=\ar
(\alpha-1)^{-1}\beta
\int^\infty_1 u^{-\beta-1}(u-1)^{1-\alpha} \dd u
=
(\alpha-1)^{-1}\beta \int^1_0 z^{\beta+\alpha-2}(1-z)^{1-\alpha}  \dd z  \cr
 \ar=\ar
\frac{\beta\Gamma(\beta+\alpha-1)\Gamma(2-\alpha)}{(\alpha-1)\Gamma(\beta+1)},
 \eeqnn
where the variable $1+z$ is replaced by $u$
and the variable $u$ is replaced by $z^{-1}$
in the second and third equality, respectively.
This gives the desired result.
\qed

In the following, we give the proof of Proposition \ref{t2.1},
which can be found in the proof of \cite[Theorem 3.2]{IkedaWatanabe}.
More precisely, we consider the function $\phi_m$  given by \eqref{3.9}
and adapt the arguments in \cite[Section 2]{XiongYang17}.
By  It\^o's formula, we consider the composite function of $\phi_m$ and the difference between two solutions before the stopping time $\gamma_n$ defined
in \eqref{3.32} of the following (see \eqref{3.14}).
We first let $m\to\infty$ (see \eqref{3.33}), then apply Gronwall's inequality
, and finally let $n\to\infty$ to conclude the proof.

{\it Proof of Proposition \ref{t2.1}.}
For $s\ge0$ define $\bar{X}_s=X_s-\tilde{X}_s$
and $\bar{Y}_s=\tilde{Y}_s-Y_s$.
For $s\ge0$, let
 \beqnn
 \ar\ar
U_{11}(s):=a_1(X_{s-}^{p_1+1}-\tilde{X}_{s-}^{p_1+1}),\quad
U_{12}(s):=\sqrt{2a_2 X_{s-}^{p_2+2}}-\sqrt{2a_2 \tilde{X}_{s-}^{p_2+2}}, \cr
 \ar\ar
U_{13}(s,u):=1_{\{u\le a_3X_{s-}^{p_3+\alpha_1}\}}-1_{\{u\le a_3\tilde{X}_{s-}^{p_3+\alpha_1}\}},\quad
U_{21}(s):=b_1(\tilde{Y}_{s-}^{q_1+1}-Y_{s-}^{q_1+1}), \cr
\ar\ar
U_{22}(s):=\sqrt{2b_2 \tilde{Y}_{s-}^{q_2+2}}-\sqrt{2b_2 Y_{s-}^{q_2+2}},\quad
U_{23}(s,u):=1_{\{u\le b_3\tilde{Y}_{s-}^{q_3+\alpha_2}\}}-1_{\{u\le b_3Y_{s-}^{q_3+\alpha_2}\}}
 \eeqnn
and
 \beqnn
V_{1}(s):=\eta_1(X_{s-}^{\theta_1}Y_{s-}^{\kappa_1}
-\tilde{X}_{s-}^{\theta_1}\tilde{Y}_{s-}^{\kappa_1}),\quad
V_2(s):=\eta_2(\tilde{Y}_{s-}^{\theta_2}\tilde{X}_{s-}^{\kappa_2}
-Y_{s-}^{\theta_2}X_{s-}^{\kappa_2}).
 \eeqnn
For $n\ge1$ define the stopping time $\gamma_n$ by
 \beqlb\label{3.32}
\gamma_n:=\tau_{n^{-1}}^-(X)\wedge\tau_{n^{-1}}^-(\tilde{X})
\wedge\tau_{n^{-1}}^-(Y)\wedge\tau_{n^{-1}}^-(\tilde{Y})\wedge
\tau_{n}^+(X)\wedge\tau_{n}^+(\tilde{X})
\wedge\tau_{n}^+(Y)\wedge\tau_{n}^+(\tilde{Y}).
 \eeqlb
By Lemma \ref{t10.1}(i), $\tau_\infty^+(X)=\tau_\infty^+(Y)=
\tau_\infty^+(\tilde{X})=\tau_\infty^+(\tilde{Y})=\infty$ almost surely.
Then $\lim_{n\to\infty}\gamma_n=\tau$ almost surely.
It follows from \eqref{a1.1} that
 \beqlb\label{3.10}
\bar{X}_{t\wedge\gamma_n}
 \ar=\ar
\bar{X}_0
-\int_0^{t\wedge\gamma_n} [U_{11}(s)+V_1(s)]\dd s
+\int_0^{t\wedge\gamma_n}U_{12}(s)\dd B_1(s) \cr
 \ar\ar
+\int_0^{t\wedge\gamma_n}\int_0^\infty\int_0^\infty z U_{13}(s,u)\tilde{N}_1(\dd s,\dd z,\dd u)
 \eeqlb
and
 \beqlb\label{3.11}
\bar{Y}_{t\wedge\gamma_n}
 \ar=\ar
\bar{Y}_0
-\int_0^{t\wedge\gamma_n} [U_{21}(s)+V_2(s)]\dd s
+\int_0^{t\wedge\gamma_n}U_{22}(s)\dd B_2(s) \cr
 \ar\ar
+\int_0^{t\wedge\gamma_n}\int_0^\infty\int_0^\infty z U_{23}(s,u)\tilde{N}_2(\dd s,\dd z,\dd u).
 \eeqlb

For $m\ge1$ define
 \beqnn
h_m:=\exp\{-m(m+1)/2\}.
 \eeqnn
Then $m^{-1}\int_{h_m}^{h_{m-1}}x^{-1}\dd x=1$.
Let $\psi_m$ be a continuous function on $\mbb{R}$
with support in $(h_m,h_{m-1})$ satisfying
 \beqnn
0\le\psi_m(x)\le
2m^{-1}x^{-1}1_{(h_m,h_{m-1})}(x),\qquad \int_{h_m}^{h_{m-1}}\psi_m(x)\dd x=1.
 \eeqnn
For $x\in\mbb{R}$ and $m\ge1$ let
 \beqlb\label{3.9}
\phi_m(x):=1_{\{x>0\}}\int_0^{x}\dd y\int_0^y\psi_m(z)\dd z.
 \eeqlb
For $m\ge1$ and $y,z\in\mbb{R}$ put
 \beqnn
\mathcal{D}_m(y,z):=\phi_m(y+z)-\phi_m(y)-z\phi_m'(y),\qquad
\mathcal{V}_m(y,z):=\phi_m(y+z)-\phi_m(y).
 \eeqnn
Let $x^+:=x\vee0$.
By \cite[Lemma 2.1]{XiongYang17} 
for all $x\in\mbb{R}$,
 \beqlb\label{3.12}
\lim_{m\to\infty}\phi_m(x)=x^+,~
\lim_{m\to\infty}\phi_m'(x)=1_{\{x>0\}},~~
|\phi_m'(x)|\le 1~\mbox{and}~|x\phi_m''(x)| \le 2/m \text{ for } m\ge1
 \eeqlb
and
 \beqlb\label{3.13}
\mathcal{D}_m(x,z)\le (2m^{-1}z^2/x)\wedge(2|z|)\quad
\text{for all}~\, m\ge1\,\text{and}~ x,z\in\mbb{R}\,\,\text{with}\,\, x(x+z)>0.
 \eeqlb

By It\^o's formula and \eqref{3.10}--\eqref{3.11} we obtain
 \beqlb\label{3.14}
 \ar\ar
\phi_m(\bar{X}_{t\wedge\gamma_n})
+\phi_m(\bar{Y}_{t\wedge\gamma_n}) \cr
 \ar=\ar
[\phi_m(\bar{X}_0)
+\phi_m(\bar{Y}_0)]
-\int_0^{t\wedge\gamma_n}
[\phi_m'(\bar{X}_{s-})U_{11}(s)+\phi_m'(\bar{Y}_{s-})U_{21}(s)]\dd s \cr
 \ar\ar
+2^{-1}\int_0^{t\wedge\gamma_n}
[\phi_m''(\bar{X}_{s-})U_{12}(s)^2+\phi_m''(\bar{Y}_{s-})U_{22}(s)^2]\dd s \cr
 \ar\ar
+\int_0^{t\wedge\gamma_n}\dd s\int_0^\infty
\Big[\int_0^\infty \mathcal{D}_m(\bar{X}_{s-}, z U_{13}(s,u))\mu_1(\dd z)
+\int_0^\infty \mathcal{D}_m(\bar{Y}_{s-}, z U_{23}(s,u))\mu_2(\dd z)\Big]\dd u \cr
 \ar\ar
-\int_0^{t\wedge\gamma_n}
[\phi_m'(\bar{X}_{s-})V_1(s)+\phi_m'(\bar{Y}_{s-})V_2(s)]\dd s \cr
 \ar\ar
+\int_0^{t\wedge\gamma_n}
\phi_m'(\bar{X}_{s-})U_{12}(s)\dd B_1(s)
+\int_0^{t\wedge\gamma_n}\phi_m'(\bar{Y}_{s-})U_{22}(s)\dd B_2(s) \cr
 \ar\ar
+\int_0^{t\wedge\gamma_n}\int_0^\infty \int_0^\infty z
\mathcal{V}_m(\bar{X}_{s-}, z  U_{13}(s,u)) \tilde{N}_1(\dd s,\dd z,\dd u) \cr
 \ar\ar
+\int_0^{t\wedge\gamma_n}\int_0^\infty \int_0^\infty z
\mathcal{V}_m(\bar{Y}_{s-}, z U_{23}(s,u)) \tilde{N}_2(\dd s,\dd z,\dd u)
 =:
I_{0,m}+\sum_{k=1}^8I_{k,m,n} .
 \eeqlb
Since $\bar{X}_0,\bar{Y}_0\le 0$, then by the definition of $\phi_m$ in \eqref{3.9},
 \beqlb\label{3.30}
I_{0,m}=0,\qquad m\ge1.
 \eeqlb
By \eqref{3.12} and the dominated convergence again,
 \beqlb\label{3.15}
\lim_{m\to\infty}\mbf{E}\big[\phi_m(\bar{X}_{t\wedge\gamma_n})\big]
+\lim_{m\to\infty}\mbf{E}\big[\phi_m(\bar{Y}_{t\wedge\gamma_n})\big]
=\mbf{E}\big[\bar{X}_{t\wedge\gamma_n}^++\bar{Y}_{t\wedge\gamma_n}^+\big].
 \eeqlb
Similarly,
 \beqlb\label{3.16}
\lim_{m\to\infty}\mbf{E}[I_{1,m,n}(t)]=-\mbf{E}\Big[\int_0^{t\wedge\gamma_n}
[1_{\{\bar{X}_{s-}>0\}}U_{11}(s)+1_{\{\bar{Y}_{s-}>0\}}U_{21}(s)]\dd s\Big]
=:I_{1,n}(t)
 \eeqlb
and
 \beqlb\label{3.17}
\lim_{m\to\infty}\mbf{E}[I_{4,m,n}(t)]=-\mbf{E}\Big[\int_0^{t\wedge\gamma_n}
[1_{\{\bar{X}_{s-}>0\}}V_1(s)+1_{\{\bar{Y}_{s-}>0\}}V_2(s)]\dd s\Big]
=:I_{4,n}(t).
 \eeqlb
By Taylor's formula,
 \beqnn
U_{12}(s)^2\le 2a_2|X_{s-}^{p_2+2}-\tilde{X}_{s-}^{p_2+2}|
\le 2a_2(p_2+2) [X_{s-}^{p_2+1}+\tilde{X}_{s-}^{p_2+1}] \cdot|\bar{X}_{s-}|
\le
4a_2(p_2+2)n^{p_2+1}|\bar{X}_{s-}|
 \eeqnn
for all $s\le\gamma_n$.
It follows from \eqref{3.12} that
 \beqnn
\int_0^{t\wedge\gamma_n}
\phi_m''(\bar{X}_{s-})U_{12}(s)^2\dd s
\le
\int_0^{t\wedge\gamma_n}
4a_2(p_2+2)n^{p_2+1}\phi_m''(\bar{X}_{s-})|\bar{X}_{s-}|\dd s
\le 8a_2(p_2+2)n^{p_2+1}m^{-1}t.
 \eeqnn
By the same arguments, we obtain
 \beqnn
\int_0^{t\wedge\gamma_n}
\phi_m''(\bar{Y}_{s-})U_{22}({s-})^2\dd s
\le 8b_2(q_2+2)n^{q_2+1}m^{-1}t.
 \eeqnn
Thus,
 \beqlb\label{3.18}
\lim_{m\to\infty}\mbf{E}[I_{2,m,n}(t)]=0.
 \eeqlb
By \eqref{4.1a},
 \beqlb\label{3.7}
|X_s^{p_3+\alpha_1}-\tilde{X}_s^{p_3+\alpha_1}|
\le (p_3+\alpha_1)|\bar{X}_s|[X_s^{p_3+\alpha_1-1}+\tilde{X}_s^{p_3+\alpha_1-1}]
\le 2(p_3+\alpha_1)n^{p_3+\alpha_1-1},s<\gamma_n.
 \eeqlb
Note that
\begin{equation}\label{3.31}
	\begin{aligned}	
U_{13}(s,u)
 &=
1_{\{a_3\tilde{X}_{s-}^{p_3+\alpha_1}< u\le a_3X_{s-}^{p_3+\alpha_1}\}}
\quad \mbox{ if }\bar{X}_{s-}>0,\\
U_{13}(s,u)
 &=
-1_{\{a_3X_{s-}^{p_3+\alpha_1}< u\le a_3\tilde{X}_{s-}^{p_3+\alpha_1}  \}}
\quad \mbox{ if }\bar{X}_{s-}<0.
\end{aligned}	
\end{equation}
It follows that $\bar{X}_{s-} U_{13}(s,u)\ge0$ for all $s,u\ge0$.
By \eqref{3.13}, for $\bar{X}_s\neq0$, we have
 \beqnn
0\le\mathcal{D}_m(\bar{X}_{s-}, z U_{13}(s,u))
\le (2m^{-1}z^2U_{13}(s,u)^2/|\bar{X}_{s-}|)\wedge(2z|U_{13}(s,u)|)
 \eeqnn
and then by \eqref{3.7} and \eqref{3.31} again, for $s<\gamma_n$,
 \beqnn
\int_0^\infty \mathcal{D}_m(\bar{X}_{s-}, z U_{13}(s,u))\dd u
 \ar\le\ar
2m^{-1}a_3z^2 |X_{s-}^{p_3+\alpha_1}-\tilde{X}_{s-}^{p_3+\alpha_1}|/|\bar{X}_{s-}| \cr
 \ar\le\ar
2 a_3(p_3+\alpha_1)n^{p_3+\alpha_1-1}  m^{-1}z^2
 \eeqnn
and
 \beqnn
\int_0^\infty \mathcal{D}_m(\bar{X}_{s-}, z U_{13}(s,u))\dd u
\le 2a_3z|X_{s-}^{p_3+\alpha_1}-\tilde{X}_{s-}^{p_3+\alpha_1}|
\le 2a_3n^{p_3+\alpha_1}z.
 \eeqnn
Thus,
 \beqnn
\int_0^\infty \mathcal{D}_m(\bar{X}_{s-}, z U_{13}(s,u))\dd u
\le2a_3(p_3+\alpha_1)n^{p_3+\alpha_1}[(m^{-1}z^2)\wedge z],\qquad s\le\gamma_n.
 \eeqnn
Similarly,
 \beqnn
\int_0^\infty \mathcal{D}_m(\bar{Y}_{s-}, z U_{23}(s,u))\dd u
\le2b_3(q_3+\alpha_2)n^{q_3+\alpha_2}[(m^{-1}z^2)\wedge z],\qquad s\le\gamma_n.
 \eeqnn
By the dominated convergence,
 \beqlb\label{3.8}
 \ar\ar
\lim_{m\to\infty}\mbf{E}[I_{3,m,n}(t)] \cr
 \ar\le\ar
\lim_{m\to\infty}2a_3(p_3+\alpha_1)n^{p_3+\alpha_1}
\mbf{E}\Big[\int_0^{t\wedge\gamma_n}\dd s\int_0^\infty
[(m^{-1}z^2)\wedge z]\mu_1(\dd z)\Big] \cr
 \ar\ar
+\lim_{m\to\infty}2b_3(q_3+\alpha_2)n^{q_3+\alpha_2}
\mbf{E}\Big[\int_0^{t\wedge\gamma_n}\dd s\int_0^\infty[(m^{-1}z^2)\wedge z]\mu_2(\dd z)\Big]=0.
 \eeqlb

By \eqref{4.1a} and \eqref{3.12},
$|\mathcal{V}_m(y,z)|\le z$ for all $y,z\in\mbb{R}$.
It follows that for $s< \gamma_n$ and $z>0$,
 \beqnn
\int_0^\infty |\mathcal{V}_m(\bar{X}_{s-}, z U_{13}(s,u))|^2\dd u
\le
z^2\int_0^\infty |U_{13}(s,u))|^2\dd u
\le
a_3z^2[X_{s-}^{p_3+\alpha_1}+\tilde{X}_{s-}^{p_3+\alpha_1}]
\le 2a_3z^2n^{p_3+\alpha_1}
 \eeqnn
and
 \beqnn
\int_0^\infty |\mathcal{V}_m(\bar{X}_{s-}, z  U_{13}(s,u))|\dd u
\le
z\int_0^\infty |U_{13}(s,u))|\dd u
\le
a_3z[X_{s-}^{p_3+\alpha_1}+\tilde{X}_{s-}^{p_3+\alpha_1}]
\le 2a_3zn^{p_3+\alpha_1}.
 \eeqnn
Thus, $t\mapsto I_{7,m,n}(t)$ is a martingale.
Similarly, $t\mapsto I_{8,m,n}(t)$ is also a martingale.
Since
 \beqnn
|\phi_m'(\bar{Y}_{s-})U_{12}(s)|^2\le
2a_2 [X_{s-}^{p_2+2}+\tilde{X}_{s-}^{p_2+2}]
\le4a_2 n^{p_2+2},\qquad s\le\gamma_n,
 \eeqnn
then by \eqref{3.12}
 \beqnn
\int_0^{t\wedge\gamma_n}|\phi_m'(\bar{Y}_{s-})U_{12}(s)|^2\dd s
\le 4a_2 n^{p_2+2}t<\infty.
 \eeqnn
Thus, $t\mapsto I_{5,m,n}(t)$ is a martingale. By the same arguments,
$t\mapsto I_{6,m,n}(t)$ is also a martingale.
Therefore,
 \beqlb\label{3.6}
\mbf{E}[I_{5,m,n}(t)+I_{6,m,n}(t)+I_{7,m,n}(t)+I_{8,m,n}(t)]=0.
 \eeqlb

Combining \eqref{3.14} with \eqref{3.30}--\eqref{3.18}
and \eqref{3.8}--\eqref{3.6} we get
 \beqlb\label{3.19}
\mbf{E}\big[\bar{X}_{t\wedge\gamma_n}^++\bar{Y}_{t\wedge\gamma_n}^+\big]
=I_{1,n}(t)+I_{4,n}(t).
 \eeqlb
Observe that
 \beqnn
1_{\{\bar{X}_{s-}>0\}}U_{11}(s)\ge0, \quad
1_{\{\bar{Y}_{s-}>0\}}U_{21}(s)\ge0,\qquad \quad s\ge0.
 \eeqnn
Thus, by \eqref{3.16},
 \beqlb\label{3.20}
I_{1,n}(t)\le0,\qquad t\ge0,~n\ge1.
 \eeqlb
If $\bar{X}_{s-}>0$ and $\bar{Y}_{s-}>0$, then
using Taylor's formula again, for $s\le\gamma_n$ we have
 \beqnn
-1_{\{\bar{X}_{s-}>0\}}V_1(s)
 \ar=\ar
\eta_1(\tilde{X}_{s-}^{\theta_1}\tilde{Y}_{s-}^{\kappa_1}-X_{s-}^{\theta_1}Y_{s-}^{\kappa_1})
\le
\eta_1 X_{s-}^{\theta_1}[\tilde{Y}_{s-}^{\kappa_1}-Y_{s-}^{\kappa_1}] \cr
 \ar\le\ar
\eta_1\kappa_1X_{s-}^{\theta_1}[\tilde{Y}_{s-}^{\kappa_1-1}+Y_{s-}^{\kappa_1-1}]\bar{Y}_{s-}
\le
2\eta_1\kappa_1n^{\theta_1}n^{|\kappa_1-1|}\bar{Y}_{s-}
 \eeqnn
and
 \beqnn
-1_{\{\bar{Y}_{s-}>0\}}V_2(s)
=\eta_2(Y_{s-}^{\theta_2}X_{s-}^{\kappa_2}
-\tilde{Y}_{s-}^{\theta_2}\tilde{X}_{s-}^{\kappa_2})
\le
\eta_2\tilde{Y}_{s-}^{\theta_2}(X_{s-}^{\kappa_2}
-\tilde{X}_{s-}^{\kappa_2})
\le
2\eta_2\kappa_2n^{\theta_2}n^{|\kappa_2-1|}\bar{X}_{s-}.
 \eeqnn
If $\bar{X}_{s-}<0$ and $\bar{Y}_{s-}>0$, then
 \beqnn
-1_{\{\bar{X}_{s-}>0\}}V_1(s)
-1_{\{\bar{Y}_{s-}>0\}}V_2(s)=
\eta_2(Y_{s-}^{\theta_2}X_{s-}^{\kappa_2}
-\tilde{Y}_{s-}^{\theta_2}\tilde{X}_{s-}^{\kappa_2})
\le
\eta_2\tilde{X}_{s-}^{\kappa_2}
[Y_{s-}^{\theta_2}-\tilde{Y}_{s-}^{\theta_2}]<0.
 \eeqnn
Therefore,
 \beqnn
-1_{\{\bar{X}_{s-}>0\}}V_1(s)
-1_{\{\bar{Y}_{s-}>0\}}V_2(s)
\le
c_n [\bar{X}_{s-}^++ \bar{Y}_{s-}^+]
 \eeqnn
with $c_n:=2[\eta_1\kappa_1n^{\theta_1+|\kappa_1-1|}
+\eta_2\kappa_2n^{\theta_2+|\kappa_2-1|}]$.
Thus, by \eqref{3.17},
 \beqnn
I_{4,n}(t)\le c_n\int_0^t
\mbf{E} \big[\bar{X}_{s\wedge\gamma_n}^++\bar{Y}_{s\wedge\gamma_n}^+ \big]\dd s.
 \eeqnn
Combining this inequality with \eqref{3.19}--\eqref{3.20} we obtain
 \beqlb\label{3.33}
\mbf{E}\big[\bar{X}_{t\wedge\gamma_n}^++\bar{Y}_{t\wedge\gamma_n}^+\big]
\le
c_n\int_0^t
\mbf{E}\big[\bar{X}_{s\wedge\gamma_n}^++\bar{Y}_{s\wedge\gamma_n}^+\big]\dd s,
\qquad t\ge0,~ n\ge1.
 \eeqlb
Now, by Gronwall's inequality,
$\bar{X}_{t\wedge\gamma_n}^+=\bar{Y}_{t\wedge\gamma_n}^+=0$ almost surely
for all $t\ge0$ and $n\ge1$.
Thus, letting $n\to\infty$ and using the
right continuity of $t\mapsto (X_t,Y_t)$ and $t\mapsto (\tilde{X}_t,\tilde{Y}_t)$
one gets the conclusion.
\qed

\noindent
{\bf Acknowledgement} The authors thank two anonymous referees for very helpful and detailed comments and suggestions.

\end{document}